\documentclass[10pt]{article}
\usepackage[latin1]{inputenc}
\usepackage{epsfig}
\usepackage{color}
\usepackage[british,english]{babel}
\usepackage{amsthm}
\usepackage{amsmath}
\usepackage{amsfonts}
\usepackage{amssymb}
\usepackage{graphicx}

\usepackage{fancyhdr}

\newtheorem{corollary}{Corollary}[section]

\newtheorem{lemma}[corollary]{Lemma}
\newtheorem{proposition}[corollary]{Proposition}
\newtheorem{remark}[corollary]{Remark}
\newtheorem{theorem}[corollary]{Theorem}

\date{}
\begin{document}
\title{On $p$-Laplacian reaction-diffusion problems with dynamical boundary conditions in perforated media}\maketitle

\vskip-30pt
 \centerline{Mar\'ia ANGUIANO
% \footnote{Departamento de An\'alisis Matem\'atico. Facultad de Matem\'aticas.
%Universidad de Sevilla, 41012 Sevilla (Spain)
%anguiano@us.es}
}

 \renewcommand{\abstractname} {\bf Abstract}
\begin{abstract} 
We study the effect of the $p$-Laplacian operator in the modelling of the heat equation through a porous medium $\Lambda_\epsilon\subset \mathbb{R}^N$ ($N\ge 2$). The case of $p=2$ was recently published in (Anguiano, Mediterr. J. Math. 17, 18 (2020)). Using rigorous functional analysis techniques and the properties of Sobolev spaces, we managed to solve additional (nontrivial) difficulties which arise compared to the study for $p=2$, and we prove a convergence theorem in appropriate functional spaces.

\end{abstract}

 {\small \bf AMS classification numbers}: 35K05, 35K58.
 
 
 \section {Statement of the problem and the results}\label{S1} 
 Homogenization problems in perforated media for the $p$-Laplacian operator have been considered in the literature over the last decades. The homogenization of the equation 
 \begin{equation}\label{elliptic_problem}
 -{\rm div}\left(|\nabla v_\epsilon|^{p-2}\nabla v_\epsilon \right)=f
 \end{equation}
in a periodically perforated domain is considered by Labani and Picard in \cite{Labani} with Dirichlet boundary conditions. Such a problem is a generalization of the linear problem for Laplace's equation, which corresponds to $p=2$, and was studied by Cioranescu and Murat in \cite{CioraMurat}. Donato and Moscariello in \cite{Donato_Moscariello} study the homogenization of a class of nonlinear elliptic Neumann problems in perforated domains of $\mathbb{R}^N$. As a consequence of \cite{Donato_Moscariello}, we are able in particular to describe the homogenization of (\ref{elliptic_problem}) with Neumann boundary conditions. This result is a generalization of earlier related works, for instance, Cioranescu and Donato \cite{Ciora2} and Cioranescu and Saint Jean Paulin \cite{Cioranescu}. In \cite{Podo_Sha2} Shaposhnikova and Podol'skii study the homogenization of (\ref{elliptic_problem}) in an $\epsilon$-periodically perforated domain with a nonlinear boundary condition. In \cite{Ildefonso1} D\'iaz {\it et al.} consider (\ref{elliptic_problem}) with a nonlinear perturbed Robin-type boundary condition in an $\epsilon$ periodically perforated domain where the size of the particles is smaller than the period $\epsilon$, and the asymptotic behavior of the solution is studied as $\epsilon\to 0$. The closest articles to this one in the literature are \cite{Gomez3, Gomez2, Gomez1}, where G\'omez {\it et al.} consider the case $2<p\leq N$, and \cite{Ildefonso2}, where D\'iaz {\it et al.} study the case $p>N$. 
%{\color{blue}In the recent book \cite{IldefonsoBook}, D\'iaz {\it et al.} collect several researches by the authors on the homogenization of nonlinear reaction-diffusion problems (mainly of elliptic or parabolic type) in the so-called ``critical scale" in which an ``anomalous" (or strange) term arises in the homogenized problem.}

It has been discovered by physicists that, as far as the Allen-Cahn equation is concerned, for certain materials a dynamical interaction with the walls must be taken into account (see Fischer {\it et al.} \cite{Fischer1,Fischer2} for more details). In this sense, in the context of heat equations, dynamical boundary conditions have been rigorously derived in Gal and Shomberg \cite{Gal2} based on first and second thermodynamical principles and their physical interpretation was also given in Goldstein \cite{Goldstein}. 
%On the other hand, in \cite{Gal_p} Gal and Warma derive the physically correct dynamical boundary condition for $p$-Laplacian reaction-diffusion equations. 
We point out that these types of boundary conditions are also used for modelling various physical situations including fluid diffusion within a semi-permeable boundary (see Crank \cite{Crank}, Langer \cite{Langer}, and, March and Weaver \cite{March} for more details) or several situations when the hear flow inside the domain is subject to nonlinear heating or cooling on the boundary (see Favini {\it et al.} \cite{Favini1,Favini2} for more details).

In the previous literature there is no study for the homogenization of $p$-Laplacian parabolic models 
%associated with nonlinear dynamical boundary conditions in a periodically perforated domain, 
as we consider in this article. Such equations model nonlinear fluid diffusion through a semi-permeable membrane (see Duvaut and Lions \cite[Ch.1]{Duvaut_Lions}) or nonlinear heat flow with radiation on the boundary causing nonlinear cooling (see Friedman \cite[Ch.7, \S 5]{Friedman}).

{\bf Model problem.} The heat equation that we study in this paper is the following
\begin{equation*}\label{PDE_1}
\displaystyle \partial_t v_\epsilon-\Delta_p\,v_\epsilon+\alpha |v_\epsilon|^{p-2}v_\epsilon =-f_1(v_\epsilon)
\quad  \text{\ in }\;\Lambda_\epsilon\times(0,\bar T) ,
\end{equation*}
where $v_\epsilon=v_\epsilon(x,t)$, $x\in \Lambda_\epsilon$, $t\in (0,\bar T)$, with $\bar T>0$ and $\alpha>0$. Assume that $\Lambda_\epsilon\subset \mathbb{R}^N$ ($N\ge 2$) is a fixed bounded domain $\Lambda$ from which a set $T_\epsilon$ of holes has been removed, in particular, $\Lambda_\epsilon$ is a periodically perforated domain with holes of the same size as the period (see \cite{Anguiano} for more details on the domain). Here the diffusion is modeled by the $p$-Laplacian operator $\Delta_p v_\epsilon:={\rm div}\left(|\nabla v_\epsilon|^{p-2}\nabla v_\epsilon \right)$ with $p\in[2,N]$. On the nonlinear term $f_1$, we assume that $f_1\in{C}\left(
\mathbb{R}\right) $, such that 
\begin{equation}\label{assumption_q_1}
p\leq q_1<+\infty,  \text{ if } p=N \quad  \text{ and } \quad 2\leq q_1\leq {Np \over N-p}, \text{ if } p\in[2,N),
\end{equation}
%\begin{equation}\label{assumption_q_2}
%2\leq q_2<+\infty,  \text{ if } p=N\quad  \text{ and } \quad 2\leq q_2\leq {(N-1)p \over N-p}, \text{ if } p\in[2,N),
%\end{equation}
\begin{equation}
\eta_{1}\left\vert s\right\vert ^{q_1}-\lambda\leq f_1(s)s\leq\eta
_{2}\left\vert s\right\vert ^{q_1}+\lambda,\quad\text{for all
$s\in\mathbb{R}$,} \label{hip_1}%
\end{equation}%
%\begin{equation}
%\alpha_{1}\left\vert s\right\vert ^{q_2}-\beta\leq g(s)s\leq\alpha
%_{2}\left\vert s\right\vert ^{q_2}+\beta,\quad\text{for all
%$s\in\mathbb{R}$,} \label{hip_2}%
%\end{equation}
and
\begin{equation}
\left( f_1(s_1)-f_1(s_2)\right) \left( s_1-s_2\right) \geq-\beta\left(
s_1-s_2\right) ^{2},\quad\text{for all
$s_1,s_2\in\mathbb{R}$,} \label{hip_3}%
\end{equation}
%and
%\begin{equation}
%\left( g(s_1)-g(s_2)\right) \left( s_1-s_2\right) \geq-l\left(
%s_1-s_2\right) ^{2},\quad\text{for all
%$s_1,s_2\in\mathbb{R}$.} \label{hip_4}%
%\end{equation}
where $\eta_i>0$, $i=1,2$, $\lambda>0$, and $\beta>0$.

We consider the following dynamical boundary conditions on the boundary of the holes
\begin{equation*}\label{PDE_2}
\partial_{\nu_p}v_\epsilon +\epsilon\,\displaystyle \partial_t v_\epsilon=-\epsilon\,f_2(v_\epsilon)  \text{\ on }%
\;\partial T_\epsilon\times( 0,\bar T),
\end{equation*}
where $T_\epsilon$ is the set of all the holes of this periodic distribution contained in $\Lambda_\epsilon$ (see \cite{Anguiano} for more details). The ``normal derivative" must be understood as $\partial_{\nu_p}v_\epsilon=|\nabla v_\epsilon|^{p-2}\nabla v_\epsilon \cdot \nu$, where $\nu$ denotes the outward normal to $\partial T_\epsilon$. This boundary equation is multiplied by $\epsilon$ to compensate the growth of the surface by shrinking $\epsilon$, where the value of $v_\epsilon$ is assumed to be the trace of the function $v_\epsilon$ defined for $x\in \Lambda_\epsilon$. On the nonlinear term $f_2$, we assume that $f_2\in{C}\left(
\mathbb{R}\right) $, such that 
%\begin{equation}\label{assumption_q_1}
%p\leq q_1<+\infty,  \text{ if } p=N \quad  \text{ and } \quad 2\leq q_1\leq {Np \over N-p}, \text{ if } p\in[2,N),
%\end{equation}
\begin{equation}\label{assumption_q_2}
2\leq q_2<+\infty,  \text{ if } p=N\quad  \text{ and } \quad 2\leq q_2\leq {(N-1)p \over N-p}, \text{ if } p\in[2,N),
\end{equation}
%\begin{equation}
%\alpha_{1}\left\vert s\right\vert ^{q_1}-\beta\leq f(s)s\leq\alpha
%_{2}\left\vert s\right\vert ^{q_1}+\beta,\quad\text{for all
%$s\in\mathbb{R}$,} \label{hip_1}%
%\end{equation}%
\begin{equation}
\eta_{1}\left\vert s\right\vert ^{q_2}-\lambda\leq f_2(s)s\leq\eta
_{2}\left\vert s\right\vert ^{q_2}+\lambda,\quad\text{for all
$s\in\mathbb{R}$,} \label{hip_2}%
\end{equation}
%\begin{equation}
%\left( f(s_1)-f(s_2)\right) \left( s_1-s_2\right) \geq-l\left(
%s_1-s_2\right) ^{2},\quad\text{for all
%$s_1,s_2\in\mathbb{R}$,} \label{hip_3}%
%\end{equation}and
and
\begin{equation}
\left( f_2(s_1)-f_2(s_2)\right) \left( s_1-s_2\right) \geq-\beta\left(
s_1-s_2\right) ^{2},\quad\text{for all
$s_1,s_2\in\mathbb{R}$.} \label{hip_4}%
\end{equation}
%where $\eta_i>0$, $i=1,2$, $\lambda>0$, and $\beta>0$.

Moreover, we consider Dirichlet boundary condition on the boundary of $\Lambda$
\begin{equation*}\label{PDE_3}
v_\epsilon= 0,  \text{\ on }%
\;\partial \Lambda\times( 0,\bar T),
\end{equation*}
and the initial conditions
\begin{equation*}\label{PDE_4}
v_\epsilon(x,0)  =  v_\epsilon^{0}(x),  \text{\ for }\;x\in\Lambda_\epsilon, \quad
v_\epsilon(x,0)  =  \phi_\epsilon^{0}(x),  \text{\ for
}\;x\in\partial T_\epsilon,
\end{equation*}
where $(v_\epsilon^{0},\phi_\epsilon^{0})$ satisfies
\begin{equation}\label{hyp 0}
v_\epsilon^{0}\in L^2\left( \Lambda\right),\quad
\phi_\epsilon^{0}\in L^{2}\left( \partial T_\epsilon\right),
\end{equation}
and
\begin{equation}\label{Initial_condition}
|v_\epsilon^0|^2_{\Lambda_\epsilon}+\epsilon|\phi_\epsilon^0|^2_{\partial T_\epsilon}\leq K,
\end{equation}
where $K>0$ and $|\cdot|_{\Lambda_\epsilon}$ (respectively $|\cdot|_{\partial T_\epsilon}$) is the norm in $L^2(\Lambda_\epsilon)$ (respectively $L^2(\partial T_\epsilon)$). Notice that on $\partial T_\epsilon$ we assume that $\phi_\epsilon^{0}(x)$ is equal to the trace of $v_\epsilon^{0}(x)$.

In summary, we study in this paper the following problem
\begin{equation}
\left\{
\begin{array}
[c]{r@{\;}c@{\;}ll}%
\displaystyle\partial_t v_\epsilon-\Delta_p\,v_\epsilon+\alpha |v_\epsilon|^{p-2}v_\epsilon &=&-f_1(v_\epsilon)
\quad & \text{\ in }\;\Lambda_\epsilon\times(0,\bar T) ,\\
\partial_{\nu_p}v_\epsilon +\epsilon\,\displaystyle\partial_t v_\epsilon&=&-\epsilon\,f_2(v_\epsilon) & \text{\ on }%
\;\partial T_\epsilon\times( 0,\bar T),\\
v_\epsilon&=& 0, & \text{\ on }%
\;\partial \Lambda\times( 0,\bar T),\\
v_\epsilon(x,0) & = & v_\epsilon^{0}(x), & \text{\ for }\;x\in\Lambda_\epsilon,\\
v_\epsilon(x,0) & = & \phi_\epsilon^{0}(x), & \text{\ for
}\;x\in\partial T_\epsilon,
\end{array}
\right. \label{PDE}%
\end{equation}
under the assumptions (\ref{assumption_q_1})-(\ref{Initial_condition}).

In a recent article (see \cite{Anguiano}) we addressed the problem (\ref{PDE}) with $p=2$ (for the physical motivation of this model, see, for instance, Timofte \cite{Timofte}). More recently, in \cite{Anguiano2} we generalize this previous study with a Laplace-Beltrami correction term
%to the case of a dynamical boundary condition of {\it reactive-diffusive} type, i.e., we add to the dynamical boundary condition a Laplace-Beltrami correction term, 
and in \cite{Anguiano3} we carried out the first study on the asymptotic behavior of the solution of parabolic models in a thin porous media.

We would like to highlight that analyzing a $p$-Laplacian problem involves additional (nontrivial) difficulties compared to the study for the Laplacian. Due to the presence of $p$-Laplacian operator in the domain, the variational formulation of the $p$-Laplacian reaction-diffusion equation is different that in \cite{Anguiano}. We have to work in the space 
\begin{eqnarray*}
V_p:=\left\{ \left( u,\gamma(u)\right) :u\in W^{1,p}\left(
\Lambda_\epsilon\right), \gamma(u)=0  \text{ on } \partial \Lambda\right\},
\end{eqnarray*}
where $\gamma$ denotes the trace operator 
$$u\in W^{1,p}(\Lambda_\epsilon)\mapsto u|_{\partial\Lambda_\epsilon}\in W^{1-{1\over p},p}(\partial\Lambda_\epsilon).$$
New $W^{1,p}$-estimates are needed to deal with problem (\ref{PDE}). To prove these estimates rigorously, we use the Galerkin approximations where we have to introduce a special basis consisting of functions in the space 
$$
\mathcal{V}_{s}:=\left\{ \left( u,\gamma(u)\right) :u\in W^{s,2}\left(
\Lambda_\epsilon\right), \gamma(u)=0  \text{ on } \partial \Lambda\right\},\quad s\ge {N(p-2)\over 2p}+1,
$$
in the sense of Lions \cite[p. 161]{Lions}. Therefore, thanks to the assumption made on $s$, we have $\mathcal{V}_{s}\subset V_p$. We use the so-called energy method introduced by Tartar \cite{Tartar} and considered by many authors (see, for instance, Cioranescu and Donato \cite{Ciora2}). Finally, in order to identity the limit equation, it is necessary to use monotonicity arguments. 
In summary, we prove that the solution of problem (\ref{PDE}), properly extended to the whole $\Lambda$, converges to the unique solution of a new nonlinear problem, defined all over the domain $\Lambda$, given by a new operator and containing extra zero order terms, capturing the effect of the influence of the non-homogeneous dynamical condition imposed on the boundary of $T_\epsilon$.
%In particular, we prove the following theorem:

\begin{theorem}[Convergence Theorem]\label{Main}
Assume (\ref{assumption_q_1})-(\ref{hip_3}), (\ref{hip_2})--(\ref{hip_4}) and (\ref{Initial_condition}). On the nonlinear term $f_2$, we suppose that $f_2\in C^1(\mathbb{R})$ and its associated exponent $q_2$ satisfies  
\begin{equation}\label{assumption_q}
p\leq q_2<+\infty \  \text{ if } \ p=N \quad \text{ and } \quad 2\leq q_2\leq {(N-1) p \over N-p} \ \text{ if } \ p\in[2,N).
\end{equation}
We suppose that $(v_\epsilon, \phi_\epsilon)$ is the unique solution of (\ref{PDE}), with $(v_\epsilon^0,\phi_\epsilon^0)\in V_p$, and where $\phi_\epsilon(t)=\gamma(v_\epsilon(t))$ a.e. $t\in (0,\bar T]$. Let $\hat v_\epsilon$ be the $W^{1,p}$-extension of $v_\epsilon$ to $\Lambda\times (0,\bar T)$. Then, we have
$$\hat v_\epsilon(t) \to v(t) \quad \text{in } L^p(\Lambda),\quad \text{as }\epsilon \to 0,\quad \forall t\in[0,\bar T],$$ 
%and $(v_\epsilon^0,\phi_\epsilon^0)\in V_p$. Let $(v_\epsilon, \phi_\epsilon)$ be the unique solution of the (\ref{PDE}), where $\phi_\epsilon(t)=\gamma(v_\epsilon(t))$ a.e. $t\in (0,\bar T]$. Then, as $\epsilon\to 0$, we have
%$$\hat v_\epsilon(t) \to v(t) \quad \text{in } L^p(\Lambda),\quad \forall t\in[0,\bar T],$$ 
where 
%$\hat \cdot$ denotes the $W^{1,p}$-extension to $\Lambda\times (0,\bar T)$, 
``$\rightarrow$'' denotes the strong convergence and $v$ is the unique solution of the problem given by
\begin{equation}\label{limit_problem_1}
\displaystyle \left(\theta^*+\theta_T \right)\displaystyle\partial_t v-{\rm div}\,b\left(\nabla v \right)+ \theta^*(\alpha |v|^{p-2}v+f_1(v))+\theta_T f_2(v)=0,
  \text{\ in }\;\Lambda\times(0,\bar T),
\end{equation}
with Dirichlet boundary condition
\begin{equation}\label{limit_problem_2}
v= 0,  \text{\ on }
\;\partial \Lambda\times( 0,\bar T),
\end{equation}
and initial condition 
\begin{equation}\label{limit_problem_3}
v(x,0)  =  v_{0}(x),  \text{\ for }\;x\in\Lambda,
\end{equation}
where $\displaystyle \theta^*=|Z^*|/ |Z|$, $\displaystyle \theta_T=|\partial T|/ |Z|$, $Z$ is the representative cell in $\mathbb{R}^N$, $T$ is an open subset of $Z$, $Z^*=Z\setminus \bar T$ and $|Z|$ (respectively $|\partial T|$ and $|Z^*|$) denotes the measure of $Z$ (respectively $\partial T$ and $Z^*$).
\\
% \begin{equation}\label{limit_problem}
%\left\{
%\begin{array}{l}
%\displaystyle \left({|Z^*|\over |Z|}+{|\partial T| \over |Z|} \right)\displaystyle\frac{\partial u}{\partial t}-{\rm div}\,b\left(\nabla u \right)+ {|Z^*|\over |Z|}(\kappa |u|^{p-2}u+f(u))+{|\partial T| \over |Z|}g(u)=0,
%  \text{\ in }\;\Lambda\times(0,T) ,\\[2ex]
%u(x,0)  =  u_{0}(x),  \text{\ for }\;x\in\Lambda,\\[2ex]
%u= 0,  \text{\ on }
%\;\partial \Lambda\times( 0,T).
%\end{array}
%\right.
%\end{equation}

For any $\zeta\in \mathbb{R}^N$, if $w(y)$ is the solution of the problem
\begin{equation}\label{system_eta}
\left\{
\begin{array}{l}
\displaystyle \int_{Z^{\star}}|\nabla_y w(y)|^{p-2}\nabla_y w(y) \cdot \nabla_y \varphi(y) dy=0\quad \forall \varphi\in  \mathbb{H}_{{\rm per}}(Z^{*}),\\[2ex]
w\in \zeta \cdot y+ \mathbb{H}_{{\rm per}}(Z^{*}),
\end{array}
\right.
\end{equation}
then $b$ is defined by
\begin{equation}\label{matrix}
b(\zeta)={1\over |Z|}\int_{Z^*}|\nabla_y w(y)|^{p-2}\nabla_y w(y) dy,
\end{equation}
where $\mathbb{H}_{{\rm per}}(Z^{*})$ is the space of functions from $W^{1,p}(Z^{*})$ which have the same trace on the opposite faces of $Z$.
\end{theorem}
\begin{remark}
If the diffusion is modeled by the Laplacian operator (i.e. $p=2$), then the limit problem (\ref{limit_problem_1})-(\ref{limit_problem_3}) is the problem obtained in \cite[Theorem 6.1]{Anguiano}.
\end{remark}
%\begin{remark}
%In this setting, a typical nonlinearity in the applications is an odd degree polynomial, 
%$\sum_{j=0}^{2k+1}c_j\,s^j,$
%where $c_{2k+1}>0$. 
%\end{remark}

\begin{remark}
It seems possible to improve the regularity assumed on the functions $f_1$ and $f_2$ when they are assumed non-increasing and Hölder continuous (for more details, see the recent book \cite{IldefonsoBook} where D\'iaz {\it et al.} present this improvement for the linear diffusion case).
\end{remark}

We organize this work as follows. In the next section, we present some notations, definitions and properties of suitable spaces for the study of (\ref{PDE}). Some preliminary results are established in Section \ref{S3}, some estimates for the solution of (\ref{PDE}) are rigorously derived in Section \ref{S4}, and a convergence result is indicated in Section \ref{S5}. Finally, in Section \ref{S6} we study the limit problem and a conclusion section is established in Section \ref{S7}.

%The article is organized as follows. In Section \ref{S2}, we introduce suitable functions spaces for our considerations. We recall here some well-know facts on the Sobolev spaces $W^{s,p}$. To prove the main result, in Section \ref{S3} we prove the existence and uniqueness of solution of (\ref{PDE}), {\it a priori} estimates are established in Section \ref{S4} and some compactness results are proved in Section \ref{S5}. Finally, the proof of Theorem \ref{Main} is established in Section \ref{S6}.

\section{The Sobolev spaces $W^{s,p}(\Lambda_\epsilon)$}\label{S2}
%{\bf Notation:} We denote by $(\cdot,\cdot) _{\Lambda_\epsilon}$ (respectively, $(
%\cdot,\cdot)_{\partial T_\epsilon}$) the inner product in
%$L^{2}(\Lambda_\epsilon)$ (respectively, in $L^{2}(\partial T_\epsilon)$),
%and by $\left\vert \cdot\right\vert _{\Lambda_\epsilon}$
%(respectively, $\left\vert \cdot\right\vert
%_{\partial T_\epsilon}$) the associated norm. We also denote by $(\cdot,\cdot) _{\Lambda_\epsilon}$ the inner product in $(L^2(\Lambda_\epsilon))^N$. 
%
%If $r\ne2$, we will also denote by
%$(\cdot,\cdot) _{\Lambda_\epsilon}$ (respectively, $(
%\cdot,\cdot)_{\partial T_\epsilon}$) the duality product between
%$L^{r'}(\Lambda_\epsilon)$ and $L^{r}(\Lambda_\epsilon)$ (respectively, the duality
%product between $L^{r'}(\partial T_\epsilon)$ and
%$L^{r}(\partial T_\epsilon)$). We will denote by
%$|\cdot|_{r,\Lambda_\epsilon}$ (respectively $|\cdot|_{r,\partial T_\epsilon}$)
%the norm in $L^r(\Lambda_\epsilon)$ (respectively in $L^r(\partial T_\epsilon)$).
%
%We denote by $(\cdot,\cdot) _{\Lambda}$ the inner product in
%$L^{2}(\Lambda)$,
%and by $\left\vert \cdot\right\vert _{\Lambda}$ the associated norm. If $r\ne2$, we will also denote by
%$(\cdot,\cdot) _{\Lambda}$ the duality product between
%$L^{r'}(\Lambda)$ and $L^{r}(\Lambda)$. We will denote by
%$|\cdot|_{r,\Lambda}$
%the norm in $L^r(\Lambda)$.
%
%%%%%%%%%%%%%%%%%%%%
%{\bf The Sobolev spaces $W^{s,p}(\Lambda_\epsilon)$:} 

In this section we recall the Sobolev spaces $W^{s,p}$, which will be used in this paper (see Adams and Fournier \cite{Adams}, Br\'ezis \cite[Chapter 9]{Brezis} and Ne${\rm \check{c}}$as \cite[Chapter 2]{Necas} for more details about them).

For any positive integer $s$ and $p\ge 1$, we define the Sobolev space $W^{s,p}(\Lambda_\epsilon)$ to be the completion of $C^{s}(\overline \Lambda_\epsilon)$, with respect to the norm
$$
||u||_{s,p,\Lambda_\epsilon}=\left(\sum_{0\leq |\alpha|\leq s}|D^{\alpha}u|^p_{p,\Lambda_\epsilon} \right)^{1/p}.
$$
Observe that $W^{s,p}(\Lambda_\epsilon)$ is a Banach space.

By the Sobolev embedding Theorem (see \cite[Chapter 4, p. 99]{Adams}), we have the embedding
\begin{equation}\label{inclusion_general}
W^{s,p}(\Lambda_\epsilon)\subset W^{1,r}(\Lambda_\epsilon),
\end{equation}
where $s\ge 1$ and ${1\over p}-{s-1\over N}\leq {1\over r}\leq {1\over p}.$

In particular, for $p\ge 1$, we define the Sobolev space $W^{1,p}(\Lambda_\epsilon)$ to be the completion of $C^{1}(\overline \Lambda_\epsilon)$, with respect to the norm
$$
||u||_{p,\Lambda_\epsilon}:=\left(|u|^p_{p,\Lambda_\epsilon}+|\nabla u|^p_{p,\Lambda_\epsilon}\right)^{1/p},
$$
where $|\cdot|_{p,\Lambda_\epsilon}$ is the norm in $L^p(\Lambda_\epsilon)$. We set $H^1(\Lambda_\epsilon)=W^{1,2}(\Lambda_\epsilon)$. 

%By $||\cdot||_{p,\Lambda_\epsilon,T}$ we denote the norm in $L^p(0,T;W^{1,p}(\Lambda_\epsilon))$. If $r\ne2$, by $|\cdot |_{r,\Lambda_\epsilon,T}$ (respectively, $|\cdot|_{r,\partial T_\epsilon,T}$), we denote the norm in $L^r(0,T;L^r(\Lambda_\epsilon))$ (respectively, $L^r(0,T;L^r(\partial T_\epsilon))$).
%
%By $\left\Vert \cdot\right\Vert _{p,\Lambda}$ we denote the norm in
%$W^{1,p}\left(\Lambda\right)$, by $||\cdot||_{p,\Lambda,T}$ we denote the norm in $L^p(0,T;W^{1,p}(\Lambda))$ and if $r\ne2$, we denote by $|\cdot|_{r,\Lambda,T}$ the norm in $L^r(0,T;L^r(\Lambda))$.

We define
\begin{equation*}
p^{\star}=\left\{
\begin{array}{l}
{Np\over N-p} \quad \text{if }p<N,\\
 +\infty \quad \text{if }p=N.
 \end{array}\right.
\end{equation*}
We have the continuous embedding
\begin{equation}\label{continuous_domain}
W^{1,p}(\Lambda_\epsilon)\subset \left\{
\begin{array}{l}
L^{r}(\Lambda_\epsilon) \quad \text{if }\quad p<N\quad \text{and} \quad r=p^{\star},\\
L^{r}(\Lambda_\epsilon)\quad \text{if }\quad p=N \quad \text{and}  \quad p\leq r<p^{\star}.
\end{array}\right.
\end{equation}
By Rellich-Kondrachov Theorem (see \cite[Chapter 9, Theorem 9.16]{Brezis}), we have the compact embedding 
\begin{equation}\label{compact_domain}
W^{1,p}(\Lambda_\epsilon)\subset \left\{
\begin{array}{l}
L^{r}(\Lambda_\epsilon) \quad \text{if }\quad p<N\quad \text{and} \quad 1\leq r< p^{\star},\\
L^{r}(\Lambda_\epsilon)\quad \text{if }\quad p=N \quad \text{and}  \quad p\leq r<p^{\star}.
 \end{array}\right.
\end{equation}
In particular, we have the compact embedding 
\begin{equation}\label{compact_domain_2}
W^{1,p}(\Lambda_\epsilon)\subset L^{2}(\Lambda_\epsilon),\ \ \forall \,2\leq p\leq N.
\end{equation}

One can define a family of spaces intermediate between $L^p$ and $W^{1,p}$. More precisely for $p\ge 1$ we define the fractional order Sobolev space
$$
W^{1-{1 \over p},p}(\partial \Lambda_\epsilon):=\left\{u\in L^{p}(\partial \Lambda_\epsilon); {|u(x)-u(y)|\over |x-y|^{1-{1\over p}+{N\over p}}}\in L^p(\partial \Lambda_\epsilon \times \partial \Lambda_\epsilon) \right\},
$$
equipped with the natural norm. We set $H^{1/2}(\partial \Lambda_\epsilon)=W^{{1\over 2},2}(\partial \Lambda_\epsilon)$. These spaces play an important role in the theory of traces.

The trace operator is denoted by $\gamma$ such that $u\mapsto
u|_{\partial\Lambda_\epsilon}$. This operator belongs to $\mathcal{L}(W^{1,p}(\Lambda_\epsilon),W^{1-{1\over p},p}(\partial\Lambda_\epsilon))$. We denote by $||\gamma||$ the norm of $\gamma$ in this space.

We will use $\|\cdot\|_{p,\partial\Lambda_\epsilon}$ to denote the
norm in $W^{1-{1\over p},p}(\partial\Lambda_\epsilon),$ which is given by
$$\|\phi\|_{p,\partial\Lambda_\epsilon}=\inf\{\|u\|_{p,\Lambda_\epsilon}:\;
\gamma(u)=\phi\}.$$ 
We define
\begin{equation*}
p^{\star}_b=\left\{
\begin{array}{l}
{(N-1)p\over N-p} \quad \text{if }p<N,\\
 +\infty \quad \text{if }p=N.
 \end{array}\right.
\end{equation*}
We have the continuous embedding
\begin{equation}\label{continuous_boundary}
W^{1-{1\over p},p}(\partial\Lambda_\epsilon)\subset \left\{
\begin{array}{l}
L^{r}(\partial \Lambda_\epsilon) \quad \text{if }\quad p<N\quad \text{and} \quad r=p^{\star}_b,\\
L^{r}(\partial \Lambda_\epsilon)\quad \text{if }\quad p=N \quad \text{and}  \quad 1\leq r<p^{\star}_b.
\end{array}\right.
\end{equation}
By \cite[Chaper 2, Theorem 6.2]{Necas}, we have the compact embedding 
\begin{equation}\label{compact_boundary}
W^{1-{1\over p},p}(\partial \Lambda_\epsilon):=\gamma\left(W^{1,p}(\Lambda_\epsilon \right))\subset L^r(\partial \Lambda_\epsilon)\quad \text{if } \quad 1\leq r<p^{\star}_b.
\end{equation}
In particular, we have the compact embedding 
\begin{equation}\label{compact_boundary_2}
W^{1-{1\over p},p}(\partial \Lambda_\epsilon)\subset L^{2}(\partial\Lambda_\epsilon) \ \ \forall \, 2\leq p\leq N.
\end{equation}
{\bf Some important notations for reading the paper:}

We define a few notations. For $\Lambda_\epsilon$, we call
$$H_p:= L^{p}\left(
\Lambda_\epsilon\right) \times L_{\partial \Lambda}^{p}\left( \partial\Lambda_\epsilon\right)
\text{,} \quad \forall p\ge 2,$$
where
$$L^p_{\partial \Lambda}(\partial \Lambda_\epsilon):=\{u\in L^p(\partial \Lambda_\epsilon):u=0  \text{ on } \partial \Lambda \}\text{,} \quad \forall p\ge 2.$$
On $H_p$, we consider the norm $|(\cdot,\cdot)|_{H_p}$ given by
$$|\left(
u,\phi\right)|^p_{H_p}=|u|_{p,\Lambda_\epsilon}^p+\epsilon|\phi|^p_{p,\partial T_\epsilon},\quad(u,\phi)\in
H_p\text{,} \quad \forall p\ge 2,$$
where
$|\cdot|_{p,\Lambda_\epsilon}$
is the norm in $L^p(\Lambda_\epsilon)$ and 
$|\cdot|_{p,\partial T_\epsilon}$
is the norm in $L^p(\partial T_\epsilon)$.

%For the sake of clarity, we shall omit to write explicitly the index $r$ if $r=2$, so 
%We denote by $H$ the Hilbert space $H:=L^{2}\left(
%\Lambda_\epsilon\right) \times L_{\partial \Lambda}^{2}\left( \partial\Lambda_\epsilon\right),$
%For each $p>1$, we consider the space
In this paper, it is very important the following space
$$V_p:=\left\{ \left( u,\gamma(u)\right) :u\in W^{1,p}_{\partial \Lambda}\left(
\Lambda_\epsilon\right) \right\}\text{,} \quad \forall p\ge 2,$$ 
where
$$W^{1,p}_{\partial \Lambda}(\Lambda_\epsilon):=\{u\in W^{1,p}(\Lambda_\epsilon):\gamma(u)=0  \text{ on } \partial \Lambda \}\text{,} \quad \forall p\ge 2.$$

Observe that $V_p$ is a closed
vector subspace of $W_{\partial \Lambda}^{1,p}\left( \Lambda_\epsilon\right) \times
W_{\partial \Lambda}^{1-{1\over p},p}\left(\partial\Lambda_\epsilon\right),$ where 
$$W^{1-{1\over p},p}_{\partial \Lambda}(\partial \Lambda_\epsilon):=\{u\in W^{1-{1\over p},p}(\partial \Lambda_\epsilon):u=0  \text{ on } \partial \Lambda \}\text{,} \quad \forall p\ge 2.$$
We endow it with the norm $\|(\cdot,\cdot)\|_{V_p}$ given by
\[
\left\Vert \left( u,\gamma(u)\right) \right\Vert^p
_{V_p}=\left\Vert
u\right\Vert _{p,\Lambda_\epsilon}^p+\left\Vert \gamma%
(u)\right\Vert^p_{p,\partial T_\epsilon}, \quad\left(
u,\gamma(u)\right) \in V_p.
\]
%Moreover, $V_p$ is a reflexive and separable space. 

Let $p'$, $q'_1$ and $q'_2$ be the conjugate exponents of $p$, $q_1$ and $q_2$, respectively. Taking into account the continuous embeddings (\ref{continuous_domain}) and (\ref{continuous_boundary}), and the assumptions (\ref{assumption_q_1})-(\ref{assumption_q_2}), we have the following useful continuous inclusions
\begin{equation}\label{inclusion1}
V_p\subset W^{1,p}(\Lambda_\epsilon)\subset L^{q_1}(\Lambda_\epsilon)\subset L^2(\Lambda_\epsilon),\quad V_p\subset W^{1-{1\over p},p}(\partial \Lambda_\epsilon)\subset L^{q_2}(\partial \Lambda_\epsilon)\subset L^2(\partial \Lambda_\epsilon),
\end{equation}
and
\begin{equation}\label{inclusion2}
L^2(\Lambda_\epsilon)\subset L^{q'_1}(\Lambda_\epsilon)\subset \left(W^{1,p}(\Lambda_\epsilon)\right)'\subset V'_p,
\end{equation}
where $\left(W^{1,p}(\Lambda_\epsilon)\right)'$ and $V'_p$ denote the dual of $W^{1,p}(\Lambda_\epsilon)$ and $V_p$, respectively. Note that $\left(W^{1,p}(\Lambda_\epsilon)\right)'$ is a subspace of $W^{-1,p'}(\Lambda_\epsilon)$, where $W^{-1,p'}(\Lambda_\epsilon)$ denotes the dual of the Sobolev space $W_0^{1,p}(\Lambda_\epsilon):= \overline{\mathcal{D}(\Lambda_\epsilon)}^{W^{1,p}(\Lambda_\epsilon)}.$

Taking into account the compact embeddings (\ref{compact_domain_2}) and (\ref{compact_boundary_2}), we have the compact embedding
\begin{equation}\label{inclusion3}
V_p\subset H_2\ \ \ \forall \, 2\leq p\leq N.
\end{equation} 
Finally, for $s\ge1$, we consider the space
$$\mathcal{V}_{s}:=\left\{ \left( u,\gamma(u)\right) :u\in W^{s,2}_{\partial \Lambda}\left(
\Lambda_\epsilon\right) \right\}.$$ We note that $\mathcal{V}_{s}$ is a closed
vector subspace of $W_{\partial \Lambda}^{s,2}\left( \Lambda_\epsilon\right) \times
W_{\partial \Lambda}^{s-{1\over 2},2}\left(\partial\Lambda_\epsilon\right).$ 

Observe that by (\ref{inclusion_general}), we have that $W^{s,2}(\Lambda_\epsilon)\subset W^{1,2}(\Lambda_\epsilon)$ and by Rellich's Theorem (see \cite[Chapter 1, Theorem 1.4]{Necas}), we obtain the compact embedding 
\begin{equation*}
W^{1,2}(\Lambda_\epsilon)\subset L^{2}(\Lambda_\epsilon),
\end{equation*}
and by \cite[Chapter 1, Exercise 1.2]{Necas}, we have the compact embedding 
\begin{equation*}
W^{s-{1\over 2},2}(\partial \Lambda_\epsilon):=\gamma_0\left(W^{s,2}(\Lambda_\epsilon \right))\subset L^2(\partial \Lambda_\epsilon).\end{equation*}
Thus, we can deduce the compact embedding
\begin{equation}\label{inclusion_Galerkin}
\mathcal{V}_{s}\subset H_2.
\end{equation}

\section{Preliminary results}\label{S3}
In the latter we will need the following results:\\
%Along this paper, we shall denote by $C$ different constants which are independent of $\epsilon$. We state in this section a result on the existence and uniqueness of solution of problem (\ref{PDE}). 
%It is easy to see from (\ref{hip_1}) and (\ref{hip_2}) that there
%exists a constant $C>0$
% such that%
\begin{remark}[Additional conditions on the nonlinear terms]Observe that we have the following conditions on the nonlinear terms:
\begin{equation}\label{hipo_consecuencia}
\left\vert f_1(s)\right\vert
 \leq
K\left( 1+\left\vert s\right\vert ^{q_1-1}\right) \text{, \ \
\ }\left\vert f_2(s)\right\vert \leq K\left( 1+\left\vert
s\right\vert ^{q_2-1}\right),\quad\text{for all
$s\in\mathbb{R}$},\quad K>0.
\end{equation}
\end{remark}

%\begin{definition}\label{definition_weakSolution} A weak solution of (\ref{PDE}) is a pair of functions $(u_\epsilon,\psi_\epsilon)$, satisfying
%\begin{equation}\label{weak0}
% u_\epsilon\in
%C([0,T];L^2(\Lambda_\epsilon)),\quad \psi_\epsilon\in
%C([0,T];L_{\partial \Lambda}^2(\partial\Lambda_\epsilon)),\quad\hbox{
%for all $T>0$,}
%\end{equation}
%\begin{equation}\label{weak1}
% u_\epsilon\in L^p(0,T;W^{1,p}(\Lambda_\epsilon)),
% \quad\hbox{
%for all $T>0$,}
%\end{equation}
%\begin{equation}\label{weak2}
%\psi_\epsilon\in L^p(0,T;W_{\partial \Lambda}^{1-{1\over p},p}(\partial\Lambda_\epsilon)),\quad\hbox{ for all $T>0$,}
%\end{equation}
%\begin{equation}\label{weak3}
% \gamma_0(u_\epsilon(t))=\psi_\epsilon(t),\quad\hbox{ a.e. $t\in (0,T],$}
% \end{equation}
% \begin{equation}\label{weak4}
%\left\{
%\begin{array}{l}
% \dfrac{d}{dt}(u_\epsilon(t),v)_{\Lambda_\epsilon}+\epsilon\,\dfrac{d}{dt}(
%\psi_\epsilon(t),\gamma_{0}(v))_{\partial T_\epsilon}+(|\nabla u_\epsilon(t)|^{p-2}\nabla u_\epsilon(t),\nabla v)_{
%\Lambda_\epsilon}+\kappa(|u_\epsilon(t)|^{p-2}u_\epsilon(t),v)_{\Lambda_\epsilon}\\[2ex]
%+(f(u_\epsilon(t)),v)_{\Lambda_\epsilon} +\epsilon\,(g(\psi_\epsilon(t)),\gamma_{0}%
%(v))_{\partial T_\epsilon}
% =0\\[2ex]
% \hbox{in $\mathcal{D}'(0,T)$, for all $v\in W_{\partial \Lambda}^{1,p}(\Lambda_\epsilon)$,}
%\end{array}
%\right.
%\end{equation}
%\begin{equation}\label{weak5}
% u_\epsilon(0)=u_\epsilon^0,\quad and\quad \psi_\epsilon(0)=\psi_\epsilon^0.
%\end{equation}
%\end{definition}

%We have the following result.

\begin{theorem}[Existence and uniqueness of solution of (\ref{PDE})]
\label{Existence_solution_PDE}Assume (\ref{assumption_q_1})--(\ref{hyp 0}), then there is a unique solution
$(v_\epsilon,\phi_\epsilon)$
of the problem (\ref{PDE}) such that, for all $\bar T>0$, $v_\epsilon\in
C([0,\bar T];L^2(\Lambda_\epsilon))\cup L^p(0,\bar T;W^{1,p}(\Lambda_\epsilon))$, $\phi_\epsilon\in
C([0,\bar T];L_{\partial \Lambda}^2(\partial\Lambda_\epsilon))\cup L^p(0,\bar T;W_{\partial \Lambda}^{1-{1\over p},p}(\partial\Lambda_\epsilon))$ where $ \gamma(v_\epsilon(t))=\phi_\epsilon(t)$ a.e. $t\in (0,\bar T]$.\\
%\begin{equation}\label{weak0}
% v_\epsilon\in
%C([0,\bar T];L^2(\Lambda_\epsilon)),\quad \phi_\epsilon\in
%C([0,\bar T];L_{\partial \Lambda}^2(\partial\Lambda_\epsilon)),\quad\hbox{
%for all $\bar T>0$,}
%\end{equation}
%\begin{equation}\label{weak1}
% v_\epsilon\in L^p(0,\bar T;W^{1,p}(\Lambda_\epsilon)),
% \quad\phi_\epsilon\in L^p(0,\bar T;W_{\partial \Lambda}^{1-{1\over p},p}(\partial\Lambda_\epsilon)),\quad\hbox{
%for all $\bar T>0$,}
%\end{equation}
%%\begin{equation}\label{weak2}
%%\psi_\epsilon\in L^p(0,T;W_{\partial \Lambda}^{1-{1\over p},p}(\partial\Lambda_\epsilon)),\quad\hbox{ for all $T>0$,}
%%\end{equation}
%with
%\begin{equation}\label{weak3}
% \gamma(v_\epsilon(t))=\phi_\epsilon(t),\quad\hbox{ a.e. $t\in (0,\bar T],$}
% \end{equation}

 We have that $(v_\epsilon,\phi_\epsilon)$ satisfies
 \begin{eqnarray}\label{weak4}
&& d_t(v_\epsilon(t),w)_{\Lambda_\epsilon}+\epsilon\,d_t(
\phi_\epsilon(t),\gamma(w))_{\partial T_\epsilon}+(|\nabla v_\epsilon(t)|^{p-2}\nabla v_\epsilon(t),\nabla w)_{
\Lambda_\epsilon}+\alpha(|v_\epsilon(t)|^{p-2}v_\epsilon(t),w)_{\Lambda_\epsilon}\\
&&+(f_1(v_\epsilon(t)),w)_{\Lambda_\epsilon} +\epsilon\,(f_2(\phi_\epsilon(t)),\gamma%
(w))_{\partial T_\epsilon}
 =0,\quad \forall w\in W_{\partial \Lambda}^{1,p}(\Lambda_\epsilon),\nonumber
 \end{eqnarray}
 in $\mathcal{D}'(0,\bar T)$,
% \begin{equation}\label{weak4}
%\left\{
%\begin{array}{l}
% \dfrac{d}{dt}(u_\epsilon(t),v)_{\Lambda_\epsilon}+\epsilon\,\dfrac{d}{dt}(
%\psi_\epsilon(t),\gamma_{0}(v))_{\partial T_\epsilon}+(|\nabla u_\epsilon(t)|^{p-2}\nabla u_\epsilon(t),\nabla v)_{
%\Lambda_\epsilon}+\kappa(|u_\epsilon(t)|^{p-2}u_\epsilon(t),v)_{\Lambda_\epsilon}\\[2ex]
%+(f(u_\epsilon(t)),v)_{\Lambda_\epsilon} +\epsilon\,(g(\psi_\epsilon(t)),\gamma_{0}%
%(v))_{\partial T_\epsilon}
% =0\\[2ex]
% \hbox{in $\mathcal{D}'(0,\bar T)$, for all $v\in W_{\partial \Lambda}^{1,p}(\Lambda_\epsilon)$,}
%\end{array}
%\right.
%\end{equation}
with the initial condition
\begin{equation}\label{weak5}
 v_\epsilon(0)=v_\epsilon^0,\quad and\quad \phi_\epsilon(0)=\phi_\epsilon^0,
\end{equation}
where $(\cdot,\cdot) _{\Lambda_\epsilon}$ is the inner product in
$L^{2}(\Lambda_\epsilon)$ or $(L^2(\Lambda_\epsilon))^N$ and the duality product between
$L^{r'}(\Lambda_\epsilon)$ and $L^{r}(\Lambda_\epsilon)$ if $r\ne 2$, and $(
\cdot,\cdot)_{\partial T_\epsilon}$ is the inner product in
$L^{2}(\partial T_\epsilon)$ and the duality product between $L^{r'}(\partial T_\epsilon)$ and
$L^{r}(\partial T_\epsilon)$ if $r\ne 2$. 

Moreover $(v_\epsilon,\phi_\epsilon)$ satisfies the energy equality
\begin{eqnarray}\nonumber
&&\frac{1}{2}d_t\left(|(v_\epsilon(t),\phi_\epsilon(t))|^2_{H_2}\right)+|\nabla
v_\epsilon(t)|^p_{p,\Lambda_\epsilon} +\alpha|v_\epsilon(t)|^p_{p,\Lambda_\epsilon} \\
 &&+(f_1(v_\epsilon(t)),v_\epsilon(t))_{\Lambda_\epsilon}+\epsilon\,(f_2(\phi_\epsilon(t)),\phi_\epsilon(t))_{\partial T_\epsilon}
 =0,\quad\mbox{a.e. $t\in(0,\bar T).$}\label{energyequality}
\end{eqnarray}
\end{theorem}
\begin{proof}
It is based on the theory of monotonicity of Lions \cite{Lions}. We will show that 
%The proof of this result is standard. For the sake of completeness, we give a sketch
%of a proof.
%
%First, we prove that $V_p$ is densely embedded in $H$. In fact, if
%we consider $\left(w,\phi\right)\in H$ such that
%\[
%(v,w)_{\Lambda_\epsilon}+\epsilon(\gamma_{0}(v),\phi)_{\partial T_\epsilon}=0,\quad
%\mbox{for all $v\in W_{\partial \Lambda}^{1,p}\left( \Lambda_\epsilon\right),$}%
%\]
%in particular, we have
%\[
%(v,w)_{\Lambda_\epsilon}=0,\quad\mbox{for all $v\in W_{0}^{1,p}\left(
%\Lambda_\epsilon\right),$}
%\]
%and therefore $w=0$. Consequently,
%\[
%(\gamma_{0}(v),\phi)_{\partial T_\epsilon}=0,\quad\mbox{for all $v\in
%W_{\partial \Lambda}^{1,p}\left( \Lambda_\epsilon\right)$,}%
%\]
%and then, as $W_{\partial \Lambda}^{1-{1\over p},p}\left( \partial\Lambda_\epsilon\right)
%=\gamma_{0}\left( W_{\partial \Lambda}^{1,p}\left( \Lambda_\epsilon\right) \right) $ is
%dense in $L_{\partial \Lambda}^{2}\left( \partial \Lambda_\epsilon\right),$ we have that
%$\phi=0$.
%Now, on the space $V_p$ we define 
%We define 
%the nonlinear
%monotone 
operator $B_p:V_p\rightarrow V_p^{\prime}$, given by
\begin{equation}\label{def_A1}
\langle B_p(( w,\gamma(w))) ,( u,\gamma (u))\rangle
:=(|\nabla w|^{p-2}\nabla w,\nabla u)_{\Lambda_\epsilon}+\alpha(|w|^{p-2}w,u)_{\Lambda_\epsilon},\quad \forall w,u\in W_{\partial \Lambda}^{1,p}\left( \Lambda_\epsilon\right),
\end{equation}
is coercive, that is $\displaystyle {\left\langle B_p\left( \left( w,\gamma(w)\right)
,\left( w,\gamma(w)\right) \right) \right\rangle \over \left\Vert \left( w,\gamma(w)\right) \right\Vert _{V_p}}\to +\infty$ when $\left\Vert \left( w,\gamma(w)\right) \right\Vert _{V_p}\to +\infty$.

We have
\begin{eqnarray}\label{Coercitivity}
\left\langle B_p\left( \left( w,\gamma(w)\right)
,\left( w,\gamma(w)\right) \right) \right\rangle &
\geq&\min\{1,\alpha\}
\left\Vert w\right\Vert _{p,\Lambda_\epsilon}^{p}%
\\
& = &\frac{1}{1+\|\gamma\|^p}\min\{1,\alpha\}
\left\Vert w\right\Vert
_{p, \Lambda_\epsilon }^{p} +\frac{\|\gamma\|^p}{1+\left\Vert
\gamma\right\Vert ^{p}} \min\{1,\alpha\}\left\Vert w\right\Vert _{p,\Lambda_\epsilon}^{p}\nonumber\\
& \geq&\frac{1}{1+\|\gamma\|^p}\min\{1,\alpha\}
 \left\Vert \left( w,\gamma(w)\right) \right\Vert _{V_p}%
^{p}\text{,}\quad \forall w\in W_{\partial \Lambda}^{1,p}(\Lambda_\epsilon),\nonumber
\end{eqnarray}so $B_p$ is coercive. By \cite[Ch.2,Th.1.4]{Lions}, we have that
\eqref{weak4}--\eqref{weak5} has a unique solution and satisfies the energy equality.

\end{proof}
\begin{remark}[Energy inequality]
By (\ref{energyequality}) and using (\ref{hip_1}) and (\ref{hip_2}), we obtain the energy inequality
\begin{eqnarray}\label{precont1}
&&d_t\left(|(v_\epsilon(t),\phi_\epsilon(t))|^2_{H_2}\right)+2|\nabla v_\epsilon(t)| _{p,\Lambda_\epsilon}^p+2\alpha| v_\epsilon(t)|^p_{p,\Lambda_\epsilon}+2\eta_1 |v_\epsilon (t)|^{q_1}_{q_1,\Lambda_\epsilon} \\
&&+2\eta_1\epsilon\,|\phi_\epsilon(t)|_{q_2,\partial T_\epsilon}^{q_2}\leq
 2\lambda\left(|\Lambda_\epsilon|+\epsilon\,|\partial T_\epsilon|\right), \nonumber
\end{eqnarray}
where $|\Lambda_\epsilon|$ (respectively $|\partial T_\epsilon|$) denotes the measure of $\Lambda_\epsilon$ (respectively $\partial T_\epsilon$).
\end{remark}
\begin{remark}[Estimates for the measures of $\Lambda_\epsilon$ and $\partial T_\epsilon$]
%Observe that the number of holes is given by $$N(\epsilon)={|\Lambda| \over (2\epsilon)^N}\left(1+o(1)\right),$$
%then using the change of variable 
%\begin{equation*}\label{dilatacion}
%y={x \over \epsilon},\quad d\sigma(y)=\epsilon^{-(N-1)}d\sigma(x),
%\end{equation*}
Taking into account the number of holes, we can deduce (see \cite[Section 4]{Anguiano2} for more details)
$$\displaystyle|\partial T_\epsilon|
%=N(\epsilon) |\partial F_{k,\epsilon}|
%=N(\epsilon)\epsilon^{N-1}|\partial F|
\leq {K \over \epsilon},\quad K>0.$$
%\end{equation*}
And since $|\Lambda_\epsilon|\leq |\Lambda|$, we have that 
\begin{equation}\label{trace2}
|\Lambda_\epsilon|+\epsilon\,|\partial T_\epsilon|\leq K,\quad K>0.
\end{equation}
\end{remark}

\section{Some estimates for the solution of (\ref{PDE})}\label{S4}
%In this section we obtain some energy estimates for the solution of (\ref{PDE}). By (\ref{energyequality}) and taking into account (\ref{hip_1})-(\ref{hip_2}), we have
%\begin{eqnarray}\label{precont1}
%&&\frac{d}{dt}\left(|(u_\epsilon(t),\psi_\epsilon(t))|^2_{H}\right)+2|\nabla u_\epsilon(t)| _{p,\Lambda_\epsilon}^p+2\kappa| u_\epsilon(t)|^p_{p,\Lambda_\epsilon}+2\alpha_1 |u_\epsilon (t)|^{q_1}_{q_1,\Lambda_\epsilon} \\
%&&+2\alpha_1\epsilon\,|\psi_\epsilon(t)|_{q_2,\partial T_\epsilon}^{q_2}\leq
% 2\beta\left(|\Lambda_\epsilon|+\epsilon\,|\partial T_\epsilon|\right), \nonumber
%\end{eqnarray}
%where $|\Lambda_\epsilon|$ and $|\partial T_\epsilon|$ denotes the measure of $\Lambda_\epsilon$ and $\partial T_\epsilon$, respectively.

%Observe that the number of holes is given by $$N(\epsilon)={|\Lambda| \over (2\epsilon)^N}\left(1+o(1)\right),$$
%then using the change of variable 
%\begin{equation*}\label{dilatacion}
%y={x \over \epsilon},\quad d\sigma(y)=\epsilon^{-(N-1)}d\sigma(x),
%\end{equation*}
%we can deduce
%\begin{equation*}
%|\partial T_\epsilon|=N(\epsilon) |\partial F_{k,\epsilon}|=N(\epsilon)\epsilon^{N-1}|\partial F|\leq {C \over \epsilon}.
%\end{equation*}
%And since $|\Lambda_\epsilon|\leq |\Lambda|$, we have that 
%\begin{equation}\label{trace2}
%|\Lambda_\epsilon|+\epsilon\,|\partial T_\epsilon|\leq C.
%\end{equation}
If we denote
\[
\mathcal{F}_1(s):=\int_{0}^{s}f_1(r)dr\quad\textrm{and}\quad \mathcal{F}_2(s):=\int_{0}^{s}f_2(r)dr,
\]
%Then, there exist positive constants $\widetilde{\alpha}_{1}$,
we can deduce that
%$\widetilde {\alpha}_{2},$ and $\widetilde{\beta}$ such that
\begin{equation}
\widetilde{\eta}_{1}| s|^{q_1}-\widetilde{\lambda} \leq\mathcal{F}_1(s)\leq\widetilde{\eta}_{2}|s|^{q_1}+\widetilde{\lambda}\quad\forall s\in\mathbb{R}, \label{hip_1_adicional}
\end{equation}
and
\begin{equation}
\widetilde{\eta}_{1}| s|^{q_2}-\widetilde{\lambda} \leq\mathcal{F}_2(s)\leq\widetilde{\eta}_{2}|s|^{q_2}+\widetilde{\lambda}\quad\forall s\in\mathbb{R}. \label{hip_2_adicional}
\end{equation}
with $\widetilde {\eta}_{i},\widetilde{\lambda}>0,$ $i=1,2$.
\begin{lemma}\label{estimates3}
We suppose (\ref{assumption_q_1})--(\ref{hip_4}) and (\ref{Initial_condition}). There is a constant $K$ independent of $\epsilon$, such that the solution $(v_\epsilon,\phi_\epsilon)$ of the problem (\ref{PDE}) satisfies
\begin{equation}\label{acotacion8_proof}
\left\Vert v_\epsilon\right\Vert _{p,\Lambda_\epsilon,\bar T}\leq K, \quad  \sup_{t\in [0,\bar T]}\left\Vert v_\epsilon(t)\right\Vert _{p,\Lambda_\epsilon}\leq K,
\end{equation}
for any initial condition $(v_\epsilon^0,\phi_\epsilon^0)\in V_p$, and where $||\cdot||_{p,\Lambda_\epsilon,\bar T}$ is the norm in $L^p(0,\bar T;W^{1,p}(\Lambda_\epsilon))$. 
%If $r\ne2$, by $|\cdot |_{r,\Lambda_\epsilon,\bar T}$, we denote the norm in $L^r(0,\bar T;L^r(\Lambda_\epsilon))$.
\end{lemma}
\begin{proof}
Using (\ref{trace2}) in (\ref{precont1}), we obtain
\begin{eqnarray}\label{precont1_new_first}
&&d_t\left(|(v_\epsilon(t),\phi_\epsilon(t))|^2_{H_2}\right)+2|\nabla v_\epsilon(t)| _{p,\Lambda_\epsilon}^p+2\alpha| v_\epsilon(t)|^p_{p,\Lambda_\epsilon}+2\eta_1 |v_\epsilon (t)|^{q_1}_{q_1,\Lambda_\epsilon} +2\eta_1\epsilon\,|\phi_\epsilon(t)|_{q_2,\partial T_\epsilon}^{q_2}\leq
 K.
\end{eqnarray}
From (\ref{def_A1})-(\ref{Coercitivity}), in particular, we can deduce
\begin{eqnarray}\label{precont1_new}
d_t\left(|(v_\epsilon(t),\phi_\epsilon(t))|^2_{H_2}\right)+\frac{2\min\{1,\alpha\}}{1+\|\gamma\|^p}|| (v_\epsilon(t),\gamma(v_\epsilon(t)))|| _{V_p}^p\leq K.
\end{eqnarray}
We integrate between $0$ and $t$ and taking into account (\ref{Initial_condition}), in particular, we have the first estimate in (\ref{acotacion8_proof}).

%%%%%%%%%%%%%%%%%%%%%%%%%%%%%%%

In order to obtain the second estimate in (\ref{acotacion8_proof}), we have to take the inner product in the problem (\ref{PDE}) with $v'_\epsilon$. To do this, we need that $v'_\epsilon \in L^p(0,\bar T;W_{\partial \Lambda}^{1,p}(\Lambda_\epsilon))\cap L^{q_1}(0,\bar T;L^{q_1}(\Lambda_\epsilon))$ with $\gamma(v'_\epsilon)\in L^{q_2}(0,\bar T;L_{\partial \Lambda}^{q_2}(\partial\Lambda_\epsilon))$. As we do not have it for our solution, we have to use the Galerkin method with the properties of $B_p$ given by (\ref{def_A1}).

At first, we introduce a special basis consisting of functions $(u_{j},\gamma(u_{j})) \in \mathcal{V}_{s}$ with $s\ge {N(p-2)\over 2p}+1$ in the sense of \cite[Chapter 2, Remark 1.6, p. 161]{Lions}. Therefore, thanks to the assumption made on $s$, taking into account (\ref{inclusion_general}), we have $\mathcal{V}_{s}\subset V_p$. The scalar product in $H_2$ generates on $\mathcal{V}_{s}\subset H_2$ the bilinear functional $ ((
w,\gamma(w)), ( u,\gamma(u)))_{H_2}$ which can be represented in the form 
$$ ((
w,\gamma(w)), ( u,\gamma(u)))_{H_2}=\langle L((w,\gamma(w))),(u,\gamma(u))\rangle_{\mathcal{V}_{s}},$$ 
where $L$ is a self-adjoint operator. The compact embedding (\ref{inclusion_Galerkin}) implies the compactness of the operator $L$. Hence, $L$ has a complete system of eigenvectors $\{(u_{j},\gamma(u_{j})) :j\geq1\}$. These vectors are orthonormal in $H_2$ and orthogonal in $\mathcal{V}_{s}$. Observe that $span\{(u_{j},\gamma(u_{j})):j\geq1\} $ is dense in $V_p$.

We use 
%by $$(u_{\epsilon,m}(t),\gamma_{0}(u_{\epsilon,m}(t)))=(u_{\epsilon,m}(t;0,u_\epsilon^0,\psi_\epsilon^0),\gamma_{0}(u_{\epsilon,m}(t;0,u_\epsilon^0,\psi_\epsilon^0)))$$ 
the Galerkin approximation of the solution of
%$(u_\epsilon(t;0,u_\epsilon^0,\psi_\epsilon^0),\gamma_{0}(u_\epsilon(t;0,u_\epsilon^0,\psi_\epsilon^0)))$ to 
(\ref{PDE}) 
%for each integer $m\geq1$, which is given by
given by
\begin{equation}
(v_{\epsilon,n}(t),\gamma(v_{\epsilon,n}(t)))=\sum_{j=1}^{n}\sigma_{\epsilon nj}(t)(u_{j},\gamma(u_{j})),\quad n\ge 1\label{Galerkin1}
\end{equation}
%and is the solution of
such that
\begin{eqnarray}
\nonumber &&d_t((v_{\epsilon,n}(t),\gamma(v_{\epsilon,n}(t))),(u_{j},\gamma(u_{j})))_{H_2}+\left\langle B_{p}((v_{\epsilon,n}(t),\gamma(v_{\epsilon,n}(t)))),(u_{j},\gamma(u_{j}))\right\rangle\\
&&+(f_1(v_{\epsilon,n}(t)),u_j)_{\Lambda_\epsilon}+\epsilon(  f_2(\gamma(v_{\epsilon,n}(t))),\gamma(u_{j}))_{\partial T_\epsilon}
=0,\quad j=1,\ldots ,n,\label{7}
\end{eqnarray}
%with initial data
\begin{equation}
\label{7'}
(v_{\epsilon,n}(0),\gamma(v_{\epsilon,n}(0)))=(v_{\epsilon,n}^{0},\gamma(v_{\epsilon,n}^{0})),
\end{equation}
and
\[
\sigma_{\epsilon nj}(t)=(v_{\epsilon,n}(t),u_{j})_{\Lambda_\epsilon}+( \gamma(v_{\epsilon,n}(t)),\gamma(u_{j}))_{\partial T_\epsilon}.
\]
%and $(u_{\epsilon,m}^{0},\gamma_0(u_{\epsilon,m}^{0}))\in span\{(w_j,\gamma_0(w_j)): j=1,\ldots ,m\}$ converges (when
%$m\to\infty$) to $(u_\epsilon^0,\psi_\epsilon^0)$ in a suitable sense which will be specified below.
%%%%%%%%%%%%%%%%%%%%%%

As $(v_\epsilon^0,\phi_\epsilon^0)\in V_{p}$, there is $(v_{\epsilon,n}^{0},\gamma(v_{\epsilon,n}^{0}))\in span\{(u_{j},\gamma(u_{j})):1\leq
j\leq n\} $, such that 
%the sequence $\{(u_{\epsilon,m}^{0},\gamma_{0}(u_{\epsilon,m}^{0}))\}$ converges to $(u_\epsilon^0,\psi_\epsilon^0)$ in $V_{p}$. Then, in particular we know that there exists a constant $C$ such that
\begin{equation}\label{Initial_condition2}
||(v_{\epsilon,n}^{0},\gamma(v_{\epsilon,n}^{0}))||_{V_p}\leq K, \quad K>0.
\end{equation}
%For each integer $m\geq1$, we consider the sequence $\{(u_{\epsilon,m}(t),\gamma_{0}(u_{\epsilon,m}(t)))\}$ defined by
%(\ref{Galerkin1})-(\ref{7'}) with these initial data.

We multiply by $\sigma'_{\epsilon nj}$ in (\ref{7}), we sum from $j=1$ to $n$, and we obtain
\begin{eqnarray}\label{equality_G}
&& |(v_{\epsilon,n}^{\prime}(t),\gamma(v_{\epsilon,n}^{\prime}(t)))|^2_{H_2}+\frac{1}{p}d_t(\left\langle B_{p}((v_{\epsilon,n}(t),\gamma(v_{\epsilon,n}(t)))),(v_{\epsilon,n}(t),\gamma(v_{\epsilon,n}(t)))\right\rangle)\nonumber\\
&& +(f_1(v_{\epsilon,n}(t)),v_{\epsilon,n}^{\prime}(t))_{\Lambda_\epsilon} +\epsilon(f_2(\gamma(v_{\epsilon,n}(t))),\gamma(v_{\epsilon,n}^{\prime}(t)))_{\partial T_\epsilon}=0.
\end{eqnarray}
%We observe that
%\[
%(f(u_{\epsilon,m}(t)),u_{\epsilon,m}^{\prime}(t))_{\Lambda_\epsilon}=\frac{d}{dt}\int_{\Lambda_\epsilon}\mathcal{F} (u_{\epsilon,m}(t))dx,
%\]
%and
%$$
%(g(\gamma_{0}(u_{\epsilon,m}(t))),\gamma_{0}(u_{\epsilon,m}^{\prime}(t)))_{\partial T_\epsilon}
%=\frac{d}{dt}\int_{\partial T_\epsilon}\mathcal{G}(\gamma_{0}(u_{\epsilon,m}(t))) d\sigma(x).
%$$
Now, we integrate between $0$ and $t$, using (\ref{Coercitivity}) and (\ref{trace2})-(\ref{hip_2_adicional}), we can deduce
\begin{eqnarray*}
&& \int_{0}^{t}|(v_{\epsilon,n}^{\prime}(s),\gamma(v_{\epsilon,n}^{\prime}(s)))|^2_{H_2}ds+\frac{\min\{1,\alpha\}}{1+\|\gamma\|^p}{1\over p} ||
(v_{\epsilon,n}(t), \gamma(v_{\epsilon,n}(t)))||^p_{V_p} \nonumber\\
&&
+\widetilde{\eta}_{1}\left(|v_{\epsilon,n}(t)|_{q_1,\Lambda_\epsilon}^{q_1} 
+\epsilon|\gamma(v_{\epsilon,n}(t))|_{q_2,\partial T_\epsilon}^{q_2} \right)
\leq {\max\{1,\alpha\} \over p} ||(v_{\epsilon,n}^{0},\gamma(v_{\epsilon,n}^{0}))||^p_{V_p}\\
&&
+\widetilde {\eta}_{2}\left(| v_{\epsilon,n}^{0}|^{q_1}_{q_1\Lambda_\epsilon}+ \epsilon|\gamma(v_{\epsilon,n}^{0})|_{q_2,\partial T_\epsilon}^{q_2}\right) +2\tilde\lambda K,\nonumber
\end{eqnarray*}
for all $t\in (0,\bar T)$. In order to estimate the right hand side of the last inequality we use (\ref{inclusion1}) together with (\ref{Initial_condition2}) and $\epsilon \ll1$. In particular, we can deduce
\begin{eqnarray}\label{last_estimate}
&& \int_{0}^{t}|(v_{\epsilon,n}^{\prime}(s),\gamma(v_{\epsilon,n}^{\prime}(s)))|^2_{H_2}ds+\frac{\min\{1,\alpha\}}{1+\|\gamma\|^p}{1\over p} ||
(v_{\epsilon,n}(t), \gamma(v_{\epsilon,n}(t)))||^p_{V_p} \leq K,
\end{eqnarray}
for all $t\in (0,\bar T)$. 
%We have proved that the sequence $\{(u_{\epsilon,m},\gamma_0(u_{\epsilon,m}))\}$ is bounded in $C([0,T];V_p),$  and $\{(u_{\epsilon,m}',\gamma_0(u_{\epsilon,m}'))\}$ is bounded in $L^2(0, T;H),$
%for all $T>0$.

%If we work with the truncated Galerkin equations (\ref{Galerkin1})-(\ref{7'}) instead of the full PDE, we note that the calculations of the proof of  (\ref{precont1_new}) can be following identically to show that $\{(u_{\epsilon,m},\gamma_0(u_{\epsilon,m}))\}$ is bounded in $L^p(0,T;V_p),$
%for all $T>0$.

Now, if we argue as in the proof of \cite[Lemma 4.5]{Anguiano}, we can deduce
%Taking into account Aubin-Lions compactness lemma (e.g., cf. Lions \cite{Lions}), the fact that $B_p$ is a monotone operator and the uniqueness of solution to (\ref{PDE}), we have that the sequence $\{(v_{\epsilon,n},\gamma(v_{\epsilon,n}))\}$ converges weakly in $
%L^{p}(0,\bar T;V_p)$ to the solution
%$(v_\epsilon,\gamma(v_\epsilon))$ to (\ref{PDE}). By (\ref{inclusion3}) and $(v_\epsilon,\gamma(v_\epsilon))\in C([0,\bar T];H_2)$, we can deduce that  
$$
\sup_{t\in [0,\bar T]}\left\Vert (v_\epsilon(t),\gamma(v_\epsilon (t)))\right\Vert _{V_p}\leq K,
$$
and, in particular, the second estimate in (\ref{acotacion8_proof}) is proved.
\end{proof}
%\section{Extension of $u_\epsilon$ to the whole $\Lambda\times (0,T)$}
%{\bf The extension of $u_\epsilon$ to the whole $\Lambda\times (0,T)$}: since the solution $u_\epsilon$ of the problem (\ref{PDE}) is defined only in $\Lambda_\epsilon\times (0,T)$, we need to extend it to the whole $\Lambda\times (0,T)$ to be able to state the convergence result. In order to do that, 
Now, we use the following extension result given by Donato and Moscariello \cite[Lemma 2.4]{Donato_Moscariello}: \\

Let $\hat v_\epsilon\in W_0^{1,p}(\Lambda)$ be a $W^{1,p}$-extension of $v_\epsilon$, that satisfies the following condition
\begin{equation}\label{Donato_extension}
|\nabla\hat v_\epsilon|_{p,\Lambda}\leq K |\nabla v_\epsilon|_{p,\Lambda_\epsilon},\quad K>0.
\end{equation}

%Using (\ref{Donato_extension}) and taking into account Lemma \ref{estimates3}, we obtain some {\it a priori} estimates for the extension of $u_\epsilon$ to the whole $\Lambda \times (0,T)$. 
\begin{corollary}\label{estimates_extension}
We suppose the assumptions in Lemma \ref{estimates3}. Then, there are a constant $K$ independend on $\epsilon$ and a $W^{1,p}$-extension $\hat v_\epsilon$ of the solution $v_\epsilon$ of (\ref{PDE}) into $\Lambda\times (0,\bar T)$, such that
\begin{equation}\label{acotacion1_extension}
\left\Vert
\hat v_\epsilon\right\Vert _{p,\Lambda,\bar T}\leq K, \quad \sup_{t\in [0,\bar T]}\left\Vert \hat v_\epsilon(t)\right\Vert _{p,\Lambda}\leq K,
\end{equation}
where $||\cdot||_{p,\Lambda,\bar T}$ is the norm in $L^p(0,\bar T;W^{1,p}(\Lambda))$ and $||\cdot||_{p,\Lambda}$ is the norm in $W^{1,p}(\Lambda)$.
%where the constant $C$ does not depend on $\epsilon$.
\end{corollary}
%\begin{proof}
%Using Poincar\'e's inequality, (\ref{Donato_extension}) and the first estimate in (\ref{acotacion8_proof}), we obtain 
%$$
%\left\Vert\tilde u_\epsilon\right\Vert^p _{p,\Lambda,T}=|\tilde u_\epsilon|^p_{p,\Lambda,T}+|\nabla \tilde u_\epsilon|^p_{p,\Lambda,T}\leq 
%C|\nabla \tilde u_\epsilon|^p_{p,\Lambda,T}\leq C|\nabla  u_\epsilon|^p_{p,\Lambda_\epsilon,T}\leq C,
%$$
%and the first estimate in (\ref{acotacion1_extension}) is proved. Similarly, using Poincar\'e 's inequality, (\ref{Donato_extension}) and the second estimate in (\ref{acotacion8_proof}), we have the second estimate in (\ref{acotacion1_extension}).
%\end{proof}

 \section{Convergence result}\label{S5}
% In this section, we obtain some compactness results about the behavior of the sequence $\tilde u_\epsilon$ satisfying the {\it a priori} estimates given in Corollary \ref{estimates_extension}.

\begin{proposition}\label{Propo_convergence}
We suppose the assumptions in Lemma \ref{estimates3}. Then, there is a function $v\in L^p(0,\bar T;W_0^{1,p}(\Lambda))$ such that, at least after extraction of a subsequence, we obtain
\begin{eqnarray}
\label{continuity1} &\hat{v}_\epsilon(t)\rightharpoonup
v(t) &\textrm{in }\ W_0^{1,p}(\Lambda),\quad \forall t\in[0,\bar T],\quad  \forall \bar T>0,
\\
\label{converge_initial_data}
&\hat{v}_{\epsilon}(t)\rightarrow
v(t)&\quad \text{in }L^p(\Lambda),\quad \forall t\in[0,\bar T],\quad  \forall \bar T>0,
\\
\label{converge_new_pfunction}
&|\hat v_{\epsilon}(t)|^{p-2}\hat v_{\epsilon}(t)\rightarrow |v(t)|^{p-2}v(t)& \quad \text{in} \quad L^{p'}(\Lambda),\quad \forall t\in[0,\bar T],\quad  \forall \bar T>0,
\\
\label{converge_f_ae}
&f_1(\hat v_{\epsilon}(t))\rightarrow f_1(v(t))& \quad \text{in} \quad L^{q'_1}(\Lambda),\quad \forall t\in[0,\bar T],\quad  \forall \bar T>0,
\end{eqnarray}
where ``$\rightharpoonup$'' denotes the weak convergence and ``$\rightarrow$'' denotes the strong convergence.\\

Moreover, there is a function $\rho\in L^{p'}(0,\bar T;L^{p'}(\Lambda))$ such that, at least after extraction of a subsequence, we obtain
\begin{eqnarray}
\label{converge_gradiente}
&\hat \rho_\epsilon \rightharpoonup
\rho &\quad \textrm{in} \quad L^{p'}(0,\bar T;L^{p'}(\Lambda)),\quad  \forall \bar T>0,
\end{eqnarray}
where $\hat \rho_\epsilon$ is given by
\begin{equation}\label{definition_tildexi}
\hat \rho_\epsilon=\left\{
\begin{array}{l}
\rho_\epsilon \quad \text{in }\Lambda_\epsilon \times (0,\bar T),\\
 0 \quad \text{in }(\Lambda\setminus \overline{\Lambda_\epsilon})\times (0,\bar T),
 \end{array}\right.
\end{equation}
with $\rho_\epsilon:=|\nabla v_\epsilon|^{p-2}\nabla v_\epsilon$.

%Let $q_2$ be the exponent satisfying (\ref{assumption_q}). Let $\bar r>1$ given by
%$$
%\bar r \in (1,p) \  \ \text{ if }p=N,\quad \bar r={Np\over (N-p)(q_2-2)+N}\ \ \text{ if } p<N.
%$$
Finally, we suppose that $f_2\in C^1(\mathbb{R})$. Then, we obtain
\begin{eqnarray}
\label{converge_g_ae}
&f_2(\hat v_{\epsilon}(t))\rightarrow f_2(v(t))& \quad \text{in} \quad L^{\bar r}(\Lambda),\quad \forall t\in[0,\bar T],\quad \forall \bar T>0,
\\
\label{converge_g2_ae}
&f_2(\hat v_{\epsilon}(t))\rightharpoonup f_2(v(t))& \quad \text{in} \quad W_0^{1,\bar r}(\Lambda),\quad \forall t\in[0,\bar T],\quad \forall \bar T>0,
\end{eqnarray}
\end{proposition}
where $\bar r>1$ is given by
$$
\bar r \in (1,p) \  \ \text{ if }p=N,\quad \bar r={Np\over (N-p)(q_2-2)+N}\ \ \text{ if } p<N,
$$
with $q_2$ satisfying (\ref{assumption_q}).
%where ``$\rightharpoonup$'' denotes weak convergence, and ``$\rightarrow$'' denotes strong convergence.
\begin{proof}
%%%%%%
If we argue as in the proof of \cite[Proposition 5.1]{Anguiano}, we obtain (\ref{continuity1})-(\ref{converge_initial_data}) and (\ref{converge_f_ae})-(\ref{converge_gradiente}).

%By the first estimate in (\ref{acotacion1_extension}), we have that
%$\{\hat v_\epsilon\}$ is bounded in $L^p(0,\bar T;W_0^{1,p}(\Lambda))$, so there exists $\{\hat v_{\epsilon'}\}\subset \{\hat v_\epsilon\}$ and  $v\in L^p(0,\bar T;W_0^{1,p}(\Lambda))$ such that
%\begin{eqnarray}
%\label{continuity1_Rellich} &\hat{v}_{\epsilon'}\rightharpoonup
%v &\textrm{in}\ L^p(0,\bar T;W_0^{1,p}(\Lambda)).
%\end{eqnarray}
%By the second estimate in (\ref{acotacion1_extension}), for each $t\in [0,\bar T]$, we get that $\{\hat v_\epsilon(t)\}$ is bounded in $W_0^{1,p}(\Lambda)$, and by (\ref{continuity1_Rellich}), we obtain (\ref{continuity1}). The conclusions of the Rellich-Kondrachov Theorem are also valid substituting $W^{1,p}(\Lambda_\epsilon)$ for $W_0^{1,p}(\Lambda)$, so using (\ref{compact_domain}) and (\ref{continuity1}), we obtain  
%\begin{equation}\label{convergence_aux_fuerte}
%\hat{v}_{\epsilon'}(t)\rightarrow
%v(t)\quad\text{in }L^r(\Lambda),\quad \forall t\in[0,\bar T],
%\end{equation}
%where $1\leq r< p^{\star}$ if $p<N$ and $p\leq r<p^{\star}$ if $p=N$. Observe that, in particular, we have (\ref{converge_initial_data}).

On the other hand, observe that
$$
||w|^{p-2}w|\leq C(1+|w|^{p-1}).
$$
Using \cite[Theorem 2.4]{Conca} with $G(x,w)=|w|^{p-2}w$, $t=p'$ and $r=p$, we can deduce that $w\in L^{p}(\Lambda)\mapsto |w|^{p-2}w\in L^{p'}(\Lambda)$ is continuous in the strong topologies. And, using (\ref{converge_initial_data}), we can deduce (\ref{converge_new_pfunction}).

%Thanks to (\ref{hipo_consecuencia}), applying \cite[Theorem 2.4]{Conca} for $G(x,w)=f_1(w)$, $t=q'_1$ and $r=q_1$, we have that the map $w\in L^{q_1}(\Lambda)\mapsto f_1(w)\in L^{q'_1}(\Lambda)$ is continuous in the strong topologies. Then, taking into account (\ref{convergence_aux_fuerte}) with $r=q_1$, we get (\ref{converge_f_ae}).
%
%From the first estimate in (\ref{acotacion8_proof}) and (\ref{definition_tildexi}), we have $|\hat \rho_\epsilon|_{p',\Lambda,\bar T}\leq C$, and hence, up to a subsequence, there exists $\rho\in L^{p'}(0,\bar T,L^{p'}(\Lambda))$ such that we obtain (\ref{converge_gradiente}).
%
%By the arbitrariness of $\bar T>0$, all the convergences are satisfied, as we wanted to prove.

%%%%%%%%%%%
For the nonlinear term $f_2$, we separate the cases $p<N$ and $p=N$. We argue as in the proof of \cite[Proposition1]{Anguiano2} in order to obtain (\ref{converge_g_ae})-(\ref{converge_g2_ae}).

\end{proof}

\section{Study of the limit problem}\label{S6}
%By $\chi_{\Lambda_\epsilon}$ we denote the characteristic function of the domain $\Lambda_\epsilon$. Due to the periodicity of the domain $\Lambda_\epsilon$, from \cite[Theorem 2.6]{CioraDonato} one has, for $\epsilon \to 0$, that 
%\begin{equation}\label{convergence_chi}
%\chi_{\Lambda_\epsilon}\stackrel{\tt
%*}\rightharpoonup {|Z^*|\over |Z|} \quad \textrm{weakly-star in}\
%L^\infty(\Lambda),
%\end{equation}
%where the limit is the proportion of the material in the cell $Z$.

%We multiply system (\ref{PDE}) by a test function $v\in \mathcal{D}(\Lambda)$, and integrating by parts, we have
%\begin{eqnarray*}
%\dfrac{d}{dt}\left(\int_{\Lambda}\chi_{\Lambda_\epsilon}\tilde u_\epsilon(t)vdx\right)+\epsilon\,\dfrac{d}{dt}\left(\int_{\partial T_\epsilon}
%\gamma_{0}(u_\epsilon(t))vd\sigma(x)\right)+\int_{\Lambda}\tilde \xi_\epsilon \cdot\nabla vdx\\[2ex]
%+\kappa \int_{\Lambda}\chi_{\Lambda_\epsilon}|\tilde u_\epsilon(t)|^{p-2}\tilde u_\epsilon(t) v dx
%+\int_{\Lambda}\chi_{\Lambda_\epsilon} f(\tilde u_\epsilon(t))vdx
% +\epsilon\, \int_{\partial T_\epsilon}
%g(\gamma_{0}(u_\epsilon(t)))vd\sigma(x)
%=0, 
% \end{eqnarray*}
% in $\mathcal{D}'(0,T)$.
 
Let $w\in \mathcal{D}(\Lambda)$ be a test function. We multiply (\ref{PDE}) by $w$ and integrating by parts, we can deduce, in $\mathcal{D}'(0,\bar T)$, the following variational formulation
\begin{eqnarray*}
d_t\left(\int_{\Lambda}\omega_{\Lambda_\epsilon}\hat v_\epsilon(t)wdx\right)+\epsilon\,d_t\left(\int_{\partial T_\epsilon}
\gamma(v_\epsilon(t))wd\sigma(x)\right)+\int_{\Lambda}\hat \rho_\epsilon \cdot\nabla wdx\\[2ex]
+\alpha \int_{\Lambda}\omega_{\Lambda_\epsilon}|\hat v_\epsilon(t)|^{p-2}\hat v_\epsilon(t) w dx
+\int_{\Lambda}\omega_{\Lambda_\epsilon} f_1(\hat v_\epsilon(t))wdx
 +\epsilon\, \int_{\partial T_\epsilon}
f_2(\gamma(v_\epsilon(t)))wd\sigma(x)
=0, 
 \end{eqnarray*}
where $\omega_{\Lambda_\epsilon}$ is the characteristic function of the domain $\Lambda_\epsilon$. Observe that the main difference with \cite[Theorem 6.1]{Anguiano} is the presence of the term $\displaystyle \int_{\Lambda}\omega_{\Lambda_\epsilon}|\hat v_\epsilon(t)|^{p-2}\hat v_\epsilon(t) w dx$.

We use $\vartheta\in C_c^1([0,\bar T])$ with $\vartheta(\bar T)=0$ and $\vartheta(0)\ne 0$. We multiply by $\vartheta$ and we integrate between $0$ and $\bar T$, in order to obtain 
%Multiplying by $\varphi$ and integrating between $0$ and $T$, we have
 \begin{eqnarray}\label{system1}\nonumber
-\vartheta(0)\left(\int_{\Lambda}\omega_{\Lambda_\epsilon}\hat v_\epsilon(0)wdx\right)-\int_0^{\bar T}\dfrac{d}{dt}\vartheta(t)\left(\int_{\Lambda}\omega_{\Lambda_\epsilon}\hat v_\epsilon(t)wdx\right)dt\\[2ex]\nonumber
-\epsilon\vartheta(0)\!\left(\int_{\partial T_\epsilon}\!\!\!
\gamma(v_\epsilon(0))wd\sigma(x)\!\right)
\!-\!\epsilon\!\!\int_0^{\bar T} \!\!\dfrac{d}{dt}\vartheta(t)\left(\int_{\partial T_\epsilon}\!\!\!
\gamma(v_\epsilon(t))wd\sigma(x)\!\right)dt\\[2ex]
+\int_0^{\bar T}\vartheta(t)\int_{\Lambda}\hat \rho_\epsilon\cdot \nabla wdxdt+\alpha \int_0^{\bar T}\vartheta(t)\int_{\Lambda}\omega_{\Lambda_\epsilon}|\hat v_\epsilon(t)|^{p-2}\hat v_\epsilon(t) w dxdt\\[2ex]\nonumber
+\int_0^{\bar T} \vartheta(t)\int_{\Lambda}\omega_{\Lambda_\epsilon}f_1(\hat v_\epsilon(t))wdxdt
 +\epsilon \int_0^{\bar T} \vartheta(t)\int_{\partial T_\epsilon}
f_2(\gamma(v_\epsilon(t)))wd\sigma(x)dt
 =0.\nonumber
 \end{eqnarray}
 
 We separately analyze the special term $\displaystyle \int_0^{\bar T}\vartheta(t)\int_{\Lambda}\omega_{\Lambda_\epsilon}|\hat v_\epsilon(t)|^{p-2}\hat v_\epsilon(t) w dxdt$.
 %For the sake of clarity, we split the proof in five parts. Firstly, we pass to the limit, as $\epsilon \to 0$, in (\ref{system1}) in order to get the limit equation satisfied by $u$. In the first step, we pass to the limit in the integrals on $\Lambda$, using Proposition \ref{Propo_convergence} and (\ref{convergence_chi}), in the second step we pass to the limit in the integrals on the boundary of the holes, where we use a convergence result based on a technique introduced by Vanninathan \cite{Vanni}, and in the third step we deduce the limit equation satisfied by $u$. In the fourth step we identify $\xi$ making use of the solutions of the problems (\ref{system_eta}), and finally we prove that $u$ is uniquely determined.

%{\bf Step 1}. Passing to the limit, as $\epsilon\to 0$, in the integrals on $\Lambda$:
%
Taking into account (\ref{converge_new_pfunction}), using \cite[Theorem 2.6]{CioraDonato} and Lebesgue's Dominated Convergence Theorem, we can deduce

%\begin{equation*}
%\int_{\Lambda}\chi_{\Lambda_\epsilon}\tilde u_\epsilon(t)v dx\to {|Z^*|\over |Z|}\int_{\Lambda}u(t)v dx,\quad \forall v \in \mathcal{D}(\Lambda),
%\end{equation*}
%\begin{equation*}
%\int_{\Lambda}\chi_{\Lambda_\epsilon}|\tilde u_\epsilon(t)|^{p-2}\tilde u_\epsilon(t) v dx\to {|Z^*|\over |Z|}\int_{\Lambda}|u(t)|^{p-2}u(t)v dx,\quad \forall v \in \mathcal{D}(\Lambda),
%\end{equation*}
%and 
%\begin{equation*}
%\int_{\Lambda}\chi_{\Lambda_\epsilon}f(\tilde u_\epsilon(t)) vdx\to {|Z^*|\over |Z|}\int_{\Lambda}f(u(t))vdx,\quad \forall v\in \mathcal{D}(\Lambda),
%\end{equation*}
%which integrating in time and using Lebesgue's Dominated Convergence Theorem, gives
%\begin{equation*}
%\int_0^T{d \over dt}\varphi(t)\left(\int_{\Lambda}\chi_{\Lambda_\epsilon}\tilde u_\epsilon(t) vdx\right)dt\to {|Z^*|\over |Z|}\int_0^T{d \over dt}\varphi(t)\left(\int_{\Lambda}u(t)vdx\right)dt,
%\end{equation*}
\begin{equation*}
\int_0^{\bar T}\vartheta(t)\int_{\Lambda}\omega_{\Lambda_\epsilon}|\hat v_\epsilon(t)|^{p-2}\hat v_\epsilon(t) wdxdt\to \theta^* \int_0^{\bar T}\vartheta(t)\int_{\Lambda}|v(t)|^{p-2}v(t)wdxdt, \quad \text{if }\quad \epsilon\to 0,
\end{equation*}
where $\theta^*=|Z^*|/ |Z|$.

For the other terms we reason as in the proof of \cite[Theorem 6.1]{Anguiano}. Then, as $\epsilon\to 0$ in (\ref{system1}), we have

%All the terms in (\ref{system1}) pass to the limit, as $\epsilon \to 0$, and therefore taking into account the previous steps, we get
 \begin{eqnarray}\label{system_previo}\nonumber
-\vartheta(0)\left(\theta^*+\theta_T \right)\left(\int_{\Lambda} v(0)wdx\right)
-\left(\theta^*+\theta_T \right)\int_0^{\bar T}\dfrac{d}{dt}\vartheta(t)\left(\int_{\Lambda} v(t)wdx\right)dt\\[2ex]
+\int_0^{\bar T}\vartheta(t)\int_{\Lambda}\rho\cdot\nabla wdxdt+\alpha\, \theta^*\int_0^{\bar T}\vartheta(t)\int_{\Lambda} |v(t)|^{p-2}v(t)w dxdt\\[2ex]
+\theta^*\int_0^{\bar T} \vartheta(t)\int_{\Lambda}f_1( v(t))wdxdt
 +\theta_T \int_0^{\bar T} \vartheta(t)\int_{\Lambda}
f_2(v(t))wdxdt
 =0,\nonumber
 \end{eqnarray}
 where $\theta_T=|\partial T|/ |Z|$.
 
 We observe that the function $\rho$ satisfies
 \begin{equation}\label{equation_xi}
 \left(\theta^*+\theta_T \right)\displaystyle\partial_t v\!-\!{\rm div}\rho\!+\! \theta^*\left(\alpha |v|^{p-2}v\!+\!f_1(v)\right)\!+\!\theta_T f_2(v)\!=\! 0, \quad \text{in} \quad \Lambda\times (0,\bar T).
 \end{equation}
 
%{\bf Step 4.} 
%It remains now to identify $\xi$. For the sake of completeness, we give a sketch of a proof, 
Following the proof of \cite[Theorem 3.1]{Donato_Moscariello}, we can deduce that
\begin{equation}\label{identificacion_xi}
\rho=b \left(  \nabla v \right)\quad \text{a.e. in }\Lambda\times (0,\bar T),
\end{equation}
where $b$ is defined by (\ref{matrix}). 

We observe that $v$ satisfies  (\ref{limit_problem_1}) using (\ref{equation_xi}) and (\ref{identificacion_xi}). The boundary condition (\ref{limit_problem_2}) is obviously satisfied. Moreover, taking into account (\ref{identificacion_xi}) in (\ref{system_previo}), we obtain exactly the variational formulation of the limit problem (\ref{limit_problem_1})-(\ref{limit_problem_3}), so we get the initial condition (\ref{limit_problem_3}). A weak solution $v$ of (\ref{limit_problem_1})-(\ref{limit_problem_3}) satisfies $v\in {C}([0,\bar T];L^{2}\left(  \Lambda\right)  )\cup L^{p}(0,\bar T;W^{1,p}_0\left(  \Lambda\right)  )$, for all $\bar T>0$, and
%\begin{equation}\label{weak1_limit}
%v\in {C}([0,\bar T];L^{2}\left(  \Lambda\right)  ),\quad v\in L^{p}(0,\bar T;W^{1,p}_0\left(  \Lambda\right)  ) \quad\text{for all }\bar T>0,
%\end{equation}
%\begin{equation}\label{weak2_limit}
%u\in L^{p}(0,T;W^{1,p}_0\left(  \Lambda\right)  ) ,\quad \text{for all }T>0,
%\end{equation}
\begin{eqnarray}\label{weak3_limit}
&&\displaystyle\left(\theta^*+\theta_T \right) d_t(v(t),w)_{\Lambda}+(b\left(\nabla v(t)\right),\nabla w)_{\Lambda}+\theta^*\alpha(|v(t)|^{p-2}v(t),w)_{\Lambda}\\
&&\displaystyle+\theta^*(f_1(v(t)),w)_{\Lambda}+\theta_T (f_2(v(t)),w)_{\Lambda}=0,\quad \forall w\in W_0^{1,p}(\Lambda),\nonumber
\end{eqnarray}
{in $\mathcal{D}'(0,\bar T)$
% \begin{equation}\label{weak3_limit}
%\left\{
%\begin{array}{l}
% \displaystyle\left({|Z^*|\over |Z|}+{|\partial F| \over |Z|} \right) \dfrac{d}{dt}(u(t),v)_{\Lambda}+(b\left(\nabla u(t)\right),\nabla v)_{\Lambda}+{|Z^*|\over |Z|}\kappa(|u(t)|^{p-2}u(t),v)_{\Lambda}\\[2ex]
%\displaystyle+{|Z^*|\over |Z|}(f(u(t)),v)_{\Lambda}+{|\partial F| \over |Z|}(g(u(t)),v)_{\Lambda}=0\\[2ex]
% \hbox{in $\mathcal{D}'(0,T)$, for all $v\in W_0^{1,p}(\Lambda)$,}
%\end{array}
%\right.
%\end{equation}
and with the initial condition
\begin{equation}\label{weak4_limit}
 v(0)=v_0,
\end{equation}
where $(\cdot,\cdot) _{\Lambda}$ is the inner product in
$L^{2}(\Lambda)$ or $(L^2(\Lambda))^N$ and the duality product between
$L^{r'}(\Lambda)$ and $L^{r}(\Lambda)$ if $r\ne 2$.

The existence and uniqueness of solution of (\ref{limit_problem_1})-(\ref{limit_problem_3}) is based on the theory of monotonicity of Lions \cite{Lions}. 
%Due to the structure properties of $b$ (see \cite[Lemmas 2.10-2.13]{Donato_Moscariello}), applying \cite[Chapter 2,Theorem 1.4]{Lions}, we obtain that the problem (\ref{limit_problem_1})-(\ref{limit_problem_3}) has a unique solution. 
%%%%%%%%%%%%%%%%%%%%
We give a sketch of a proof.

We take into account the structure properties of $b$, in particular, from \cite[Lemmas 2.10-2.13]{Donato_Moscariello}), for any $\zeta\in \mathbb{R}^N$, we have
\begin{equation}\label{propiedad1}
|b(\zeta)|\leq c(1+|\zeta|)^{p-1},
\end{equation}
where $c>0$, and for $\zeta_1$, $\zeta_2\in \mathbb{R}^N$, we have
\begin{equation}\label{propiedad2}
(b(\zeta_1)-b(\zeta_2),\zeta_1-\zeta_2)\ge \kappa|\zeta_1-\zeta_2|^p,\quad \text{if }p\ge2,
\end{equation}
where $\kappa>0$.

On the space $W_0^{1,p}(\Lambda)$ we define the nonlinear
monotone operator $B_p:W_0^{1,p}(\Lambda)\rightarrow (W_0^{1,p}(\Lambda))^{\prime}$, given by
\begin{equation}\label{def_A1_limit}
\langle B_p(w) ,u\rangle
:=C_1(b(\nabla w),\nabla u)_{\Lambda}+C_1\theta^*\alpha(|w|^{p-2}w,u)_{\Lambda},\quad \forall w,u\in W_{0}^{1,p}\left( \Lambda\right),
\end{equation}
where $C_{1}=\displaystyle\left(\theta^*+\theta_T \right)^{-1}$.

Taking into account (\ref{propiedad2}) with $\zeta_1=\nabla w$ and $\zeta_2=0$, we can deduce
\begin{eqnarray}\label{Coercitivity_limit}
\left\langle B_p\left( w\right)
,w  \right\rangle &
=&C_1(b(\nabla w),\nabla w)_{\Lambda}+C_1\theta^*\alpha|w|^{p}_{p,\Lambda}
\\
& \geq &C_1\kappa |\nabla w|^p_{p,\Lambda}+C_1\int_{\Lambda}b(0)\nabla w dx+C_1\theta^* \alpha|w|^{p}_{p,\Lambda}\text{,}\quad \forall w\in W_{0}^{1,p}(\Lambda),\nonumber
\end{eqnarray}
where
$|\cdot|_{p,\Lambda}$
is the norm in $L^p(\Lambda)$. 

On the other hand, using (\ref{propiedad1}) with $\zeta=0$ and Young's inequality, we can deduce
\begin{eqnarray}\label{Young}
\left|\int_{\Lambda}b(0)\nabla w dx \right| &
\leq&\int_{\Lambda}c|\nabla w|dx\leq {2\over p\kappa}{c^{p'}\over p'}|\Lambda|+{\kappa\over 2}|\nabla w|^p_{p,\Lambda},
\end{eqnarray}
where $|\Lambda|$ denotes the measure of $\Lambda$ and $p'$ is the conjugate exponent of $p$.

Then, taking into account (\ref{Young}) in (\ref{Coercitivity_limit}), we obtain
\begin{eqnarray*}
\left\langle B_p\left( w\right)
,w  \right\rangle +C_1{2\over p\kappa}{c^{p'}\over p'}|\Lambda|
 &\geq& C_1{\kappa\over 2} |\nabla w|^p_{p,\Lambda}+C_1\theta^*\alpha|w|^{p}_{p,\Lambda}\\
&\ge &\min\{C_1{\kappa\over 2},C_1\theta^*\alpha\}||w||^p_{p,\Lambda}\text{,}\quad \forall w\in W_{0}^{1,p}(\Lambda),\nonumber
\end{eqnarray*} so $B_p$ is coercive.

Now, we consider the following spaces and operators
$$ V_1=W_0^{1,p}(\Lambda),\quad V_{2}=L^{q_1}( \Lambda),\quad
V_{3}= L^{q_2}(\Lambda),$$
$$B_1(v)=B_p,\quad B_{2}(v) =C_1\theta^*f_1(w),\quad
B_{3}(v) =C_1\theta_T f_2(w).
$$

We observe that from \eqref{hipo_consecuencia}, we can deduce that $B_i:V_i\rightarrow V'_i$ for $i=2,3.$

Taking into account the continuous embedding (\ref{continuous_domain}) for $\Lambda$, and the assumptions (\ref{assumption_q_1}) and (\ref{assumption_q}), we have the following useful continuous inclusions
\begin{equation}\label{inclusion1_limit}
W_0^{1,p}(\Lambda)\subset W^{1,p}(\Lambda)\subset L^{q_1}(\Lambda)\subset L^2(\Lambda),\quad W_0^{1,p}(\Lambda)\subset W^{1,p}(\Lambda)\subset L^{q_2}(\Lambda)\subset L^2(\Lambda),
\end{equation}
and
\begin{equation}\label{inclusion2_limit}
L^2(\Lambda)\subset L^{q'_1}(\Lambda)\subset \left(W^{1,p}(\Lambda)\right)'\subset \left(W_0^{1,p}(\Lambda)\right)',\quad L^2(\Lambda)\subset L^{q'_2}(\Lambda)\subset \left(W^{1,p}(\Lambda)\right)'\subset \left(W_0^{1,p}(\Lambda)\right)',
\end{equation}
where $\left(W^{1,p}(\Lambda)\right)'$ and $\left(W_0^{1,p}(\Lambda)\right)'$ denote the dual of $W^{1,p}(\Lambda)$ and $W_0^{1,p}(\Lambda)$, respectively.

%With this notation, and denoting $V=\cap_{i=1}^3V_i=W_0^{1,p}(\Lambda),$ $p_1=p,$
%$p_2=q_1,$ $p_3=q_2,$ one has that
%\eqref{weak1_limit}--\eqref{weak4_limit} is equivalent to
%\begin{equation}\label{Weak1_limit}
% u\in C([0,T];L^2(\Lambda)),\quad u\in\bigcap_{i=1}^3 L^{p_i}(0,T;V_i)=L^{p}(0,T;W^{1,p}_0\left(  \Lambda\right)  ) ,
% \quad\mbox{for all $T>0,$}
%\end{equation}
%\begin{equation}\label{Weak2_limit}
%u'(t)+\sum_{i=1}^3A_i(u(t))=0\quad\mbox{in
%$\mathcal{D}'(0,T;V'),$}
%\end{equation}
%\begin{equation}\label{Weak3_limit}
%u(0)=u_0.
%\end{equation}

%Taking into account (\ref{inclusion1_limit})-(\ref{inclusion2_limit}) and 
Finally, if we apply \cite[Ch.2,Th.1.4]{Lions}, we have that \eqref{weak3_limit}--\eqref{weak4_limit} has a unique solution. As $v$ is uniquely determined, the whole sequence $\hat v_\epsilon$ converges to $v$ and this completes the proof of Theorem \ref{Main}.

\section{Conclusions}\label{S7}
In this paper, we consider a parabolic model in a perforated media $\Lambda_\epsilon \subset \mathbb{R}^N$ ($N\ge 2$) with periodically distributed holes of size $\epsilon$. The $p$-Laplacian operator appears in wide range of scientific fields, for instance in fluid dynamics (e.g. flow in a porous media), nonlinear elasticity, glaciology and image restoration. In this sense, in $\Lambda_\epsilon\times (0,\bar T)$, with $\bar T>0$, we consider the $p$-Laplace heat equation 
$$
\displaystyle \partial_t v_\epsilon-\Delta_p\,v_\epsilon+\alpha |v_\epsilon|^{p-2}v_\epsilon =-f_1(v_\epsilon),
$$
where 
$$
\Delta_p v_\epsilon:={\rm div}\left(|\nabla v_\epsilon|^{p-2}\nabla v_\epsilon \right),$$
 and $f_1$ is the derivative of a potential function $\mathcal{F}_1$, that is, $\mathcal{F}_1(s)=\int_0^s f_1(r)dr$.

The usual boundary conditions considered in the literature are Dirichlet or Neumann. With this standard boundary condition, if we consider the following energy functional 
\begin{equation}\label{energy1}
\mathcal{E}_{\Lambda_\epsilon}(v_\epsilon(t)):=\int_{\Lambda_\epsilon}\left({1\over p}|\nabla v_\epsilon(t)|^p+{\alpha\over p}| v_\epsilon(t)|^p+\mathcal{F}_1(v_\epsilon(t)) \right)dx,
\end{equation}
and taking into account that
\[
d_t\left(\int_{\Lambda_\epsilon}\mathcal{F}_1 (v_{\epsilon}(t))dx\right)=(f_1(v_{\epsilon}(t)),v_{\epsilon}^{\prime}(t))_{\Lambda_\epsilon},
\]
then, we observe that this energy functional is decreasing.

In this paper, instead of Dirichlet or Neumann boundary condition, we consider a boundary condition, which depends on the time, on the boundary of the holes. In this sense, we add to (\ref{energy1}) the following energy functional
\begin{equation}\label{energy2}
\mathcal{E}_{\partial T_\epsilon}(v_\epsilon(t)):=\int_{\partial T_\epsilon}\mathcal{F}_2(v_\epsilon(t))dx,
\end{equation}
where $T_\epsilon$ is the set of all the holes contained in a bounded open set $\Lambda\subset \mathbb{R}^N$ and $\mathcal{F}_2$ is a nonlinear function such that $\mathcal{F}_2=\int_0^s
f_2(r)dr$. Then, we obtain a total energy functional 
$$\mathcal{E}(v_\epsilon(t))=\mathcal{E}_{\Lambda_\epsilon}(v_\epsilon(t))+\mathcal{E}_{\partial T_\epsilon}(v_\epsilon(t)),$$ which, using integration by parts, is decreasing for all time $t\ge 0$ if we consider the following dynamic boundary condition on the boundary of the holes
$$
\partial_{\nu_p}v_\epsilon +\displaystyle \partial_t v_\epsilon=-\,f_2(v_\epsilon),
$$
where $\partial_{\nu_p}v_\epsilon=|\nabla v_\epsilon|^{p-2}\nabla v_\epsilon \cdot \nu$, with $\nu$ the outward normal to $\partial T_\epsilon$, and Dirichlet boundary condition on the boundary of $\Lambda$.

We extend the results of \cite{Anguiano} obtained for $p=2$ to the case $2\leq p \leq N$. The main result of this paper (Theorem \ref{Main}) could be summarized by the following expansion for $v_\epsilon$
$$\hat v_\epsilon \sim v,$$ 
where $\hat v_\epsilon$ is the $W^{1,p}$-extension of $v_\epsilon$ to $\Lambda$ and $v$ is the solution of a parabolic model coming the homogenization in the porous media.

Using the present study as a starting point, various improvements can be proposed. The first one is the generalization of the asymptotic study to a other types of nonlinear diffusion (and not only $p$-Laplacian operator). For instance, it is very interesting if the operators $\Delta_p v_\epsilon$ and $\partial_{\nu_p} v_\epsilon$ from (\ref{PDE}) are replaced by the operators ${\rm div}\left(a\left(|\nabla v_\epsilon| \right)\nabla v_\epsilon \right)$ and $b(x)a\left(|\nabla v_\epsilon|\right)\nabla v_\epsilon\cdot \nu$, where $b\in L^{\infty}(\partial T_\epsilon)$, $b\ge b_0>0$ and $a\in C^1\left(\mathbb{R}^N,\mathbb{R}\right)$ is a monotone nondecreasing function such that there are two positive constants $c_1$, $c_2$ such that 
$$|a(y)|\leq c_1(1+|y|^{p-2}),\quad a(y)|y|^2\ge c_2|y|^p,\quad \forall y\in \mathbb{R}^N.$$
Another possible way is to study this parabolic model in a thin porous media (see, for instance, \cite{Anguiano4, Anguiano6,Anguiano5,Anguiano7,Grau1,Grau2} for more details on the importance of this type of domains). Finally, another problem could be to consider a porous media containing a thin fissure. This type of domain is very interesting because it models cracks in geological strata (see, for instance, \cite{Anguiano10,Anguiano9,Anguiano8,Anguiano11} for more details). Mathematical models of such domains include several small parameters, one is connected to the domain height or the width of a thin fissure and the others to the microstructure. This approach could be very interesting.

\subsection*{Acknowledgments}
María Anguiano is very grateful to the anonymous referees for their support of this document with their nice and valuable comments and suggestions.

\end{document}